\documentclass[12pt]{scrartcl}

\usepackage[
  bookmarks,
  bookmarksopen=true,
  bookmarksnumbered=true,
  pdfusetitle,
  pdfcreator={},
  colorlinks,
  linkcolor=black,
  urlcolor=black,
  citecolor=black,
  plainpages=false,
  ]{hyperref}

\usepackage{amsfonts}
\usepackage{amsmath}
\usepackage{amssymb}
\usepackage{amsthm}
\usepackage{algorithm}
\usepackage{algorithmic}
\usepackage{graphicx}
\usepackage{color}
\usepackage{subcaption}
\usepackage{tikz,pgfplots}
\usepackage{prettyref}
\usepackage{mathrsfs}	
\usepackage[utf8]{inputenc}

\newcommand{\N}{\ensuremath{\mathbb{N}_0}}

\renewcommand{\S}{\ensuremath{\mathbb{S}}}
\newcommand{\B}{\ensuremath{\mathbb{B}}}

\newcommand{\R}{\ensuremath{\mathbb{R}}}
\newcommand{\C}{\ensuremath{\mathbb{C}}}

\newcommand{\abs}[1]{\ensuremath{\left\vert#1\right\vert}}
\newcommand{\dx}{\mathrm{d}}

\newcommand{\e}{\mathrm{e}}

\newcommand{\zb}[1]{\ensuremath{\boldsymbol{#1}}}

\renewcommand{\d}{\, \mathrm{d}}

\newcommand{\bigo}{\ensuremath{\mathcal{O}}}

\newcommand{\norm}[1]{\left\lVert \smash{#1} \right\rVert}

\newcommand\multiset[2]
{\mathchoice{\left(\kern-0.5em{\binom{#1}{#2}}\kern-0.5em\right)}
	{\bigl(\kern-0.3em{\binom{#1}{#2}}\kern-0.3em\bigr)}
	{\bigl(\kern-0.3em{\binom{#1}{#2}}\kern-0.3em\bigr)}
	{\bigl(\kern-0.3em{\binom{#1}{#2}}\kern-0.3em\bigr)}}

\newtheorem{theorem}{Theorem}[section]
\newtheorem{lemma}[theorem]{Lemma}
\newtheorem{remark}[theorem]{Remark}
\newtheorem{definition}[theorem]{Definition}
\newtheorem{example}[theorem]{Example}
\newtheorem{corollary}[theorem]{Corollary}
\newtheorem{proposition}[theorem]{Proposition}
\newtheorem{problem}[theorem]{Problem}

\theoremstyle{remark}
\newtheorem*{remark*}{Remark}

\newenvironment{Theorem}[1][noisnotdefined]{ \ifthenelse{\equal{#1}{noisnotdefined}}{\begin{theorem}}{\begin{theorem}[#1]}\normalfont\slshape}{\end{theorem}}
\newenvironment{Lemma}{ \begin{lemma}\normalfont\slshape}{\end{lemma}}
\newenvironment{Remark}[1][noisnotdefined]{ \ifthenelse{\equal{#1}{noisnotdefined}}{\begin{remark}}{\begin{remark}[#1]}\normalfont\rmfamily}{\bend\end{remark}}

\numberwithin{equation}{section}

\newcommand{\bend}{\hspace*{0ex} \hfill \hbox{\vrule height
    1.5ex\vbox{\hrule width 1.4ex \vskip 1.4ex\hrule  width 1.4ex}\vrule
    height 1.5ex}}

\long\def\symbolfootnote[#1]#2{\begingroup
\def\thefootnote{\fnsymbol{footnote}}\footnote[#1]{#2}\endgroup}

\newcounter{todocounter}
\newcommand{\todo}[2][noisnotdefined]{
 \marginpar{\fcolorbox{black}{yellow}{\footnotesize\textbf{todo}}
 \ifthenelse{\equal{#1}{noisnotdefined}}{}{\textcolor{black}{\newline\tiny #1}}}
 \textbf{\ifthenelse{\equal{#2}{.}}
   {\fcolorbox{red}{white}{\textcolor{red}{$\maltese$}}}{{\textcolor{red}{#2}}}}
 \refstepcounter{todocounter}}

\newrefformat{def}{Definition \ref{#1}}
\newrefformat{sub}{Subsection \ref{#1}}
\newrefformat{prop}{Proposition \ref{#1}}
\newrefformat{fig}{Figure \ref{#1}}
\newrefformat{rem}{Remark \ref{#1}}
\newrefformat{cor}{Corollary \ref{#1}}

\title{The cone-beam transform and spherical convolution operators}

\date{\today}
\date{}
\author{ Michael Quellmalz\thanks{Faculty of Mathematics, Chemnitz University of Technology, D-09107 Chemnitz, Germany.\newline E-mail:		\href{mailto:michael.quellmalz@math.tu-chemnitz.de}{michael.quellmalz@math.tu-chemnitz.de}}
	\and Ralf Hielscher\thanks{Faculty of Mathematics, Chemnitz University of Technology, D-09107 Chemnitz, Germany.\newline E-mail:		\href{mailto:ralf.hielscher@math.tu-chemnitz.de}{ralf.hielscher@math.tu-chemnitz.de}}
	\and Alfred K. Louis\thanks{Department of Mathematics, Saarland University, D-66041 Saarbrücken, Germany.\newline E-mail:		\href{mailto:louis@num.uni-sb.de}{louis@num.uni-sb.de}}
}

\begin{document}

\maketitle

\begin{abstract}

The cone-beam transform consists of integrating a function defined on the three-dimensional space along every ray that starts on a certain scanning set. Based on Grangeat's formula, Louis [2016, \emph{Inverse Problems} \textbf{32} 115005] states reconstruction formulas based on a new generalized Funk--Radon transform on the sphere.

In this article, we give a singular value decomposition of this generalized Funk--Radon transform.
We use this result to derive a singular value decomposition of the cone-beam transform with sources on the sphere thus generalizing a result of Kazantsev [2015, \emph{J.\ Inverse Ill-Posed Probl.} {23}(2):173--185].

\medskip

\end{abstract}

\section{Introduction}

The cone-beam transform integrates a function $f\colon \R^d\to\R$ along every ray that starts in a certain scanning set $\Gamma\subset\R^d$. We define the cone-beam transform, or divergent beam X-ray transform, by
\begin{equation*}
\mathcal Df (\zb a,\zb\omega)
= \int_{0}^{\infty} f(\zb a + t\zb\omega)\d t
,\qquad \zb\omega\in\S^{d-1},\; \zb a\in\Gamma.
\end{equation*}
The cone-beam transform is widely used in medical imaging and nondestructive testing of three-dimensional objects.
Uniqueness of the reconstruction was shown under rather weak assumptions, namely that $\Gamma$ is an infinite set with positive distance to the convex hull of the support of $f$, see \cite{HaSmSoWa80}.
An explicit inversion formula \cite{Tuy83} is known in case the Tuy--Kirillov completeness condition is satisfied, which states that the scanning set $\Gamma$ intersects every hyperplane hitting $\operatorname{supp} f$ transversally, see \cite[Chap.\ 2]{NaWue00}.
In 3D, if the Tuy--Kirillov condition is not satisfied, one can stably detect singularities of~$f$ only along planes that meet the scanning curve $\Gamma$, see \cite{Qui93}.

The present article's main focus is the setting where the scanning set $\Gamma$ covers the whole sphere, where, however, in most practical application, $\Gamma$ is considered to be only a curve, cf.\ \cite{Fin85,Smi90}.
In 2015, Kazantsev \cite{Kaz15} showed the singular value decomposition of the cone-beam transform $\mathcal D$ where the function $f$ is supported inside the unit ball $\B^3$ and the scanning set $\Gamma$ is the sphere~$\S^2$.
The singular value decomposition of the parallel beam X-ray transform is due to Maaß \cite{Maa87}.
In 2016, Louis \cite{Lou16} gave new inversion formulas for the cone-beam transform for dimension $d=3$.
The proof of these formulas utilized Grangeat's formula and a generalized Funk--Radon transform $\mathcal S^{(j)}$ on the sphere,
which is defined for a function $f\colon\S^{d-1}\to\C$ by
\begin{equation*}
\mathcal S^{(j)}_d f (\zb \xi)
= \int_{\S^{d-1}} \delta^{(j)}(\zb\xi^\top\zb\eta)\, f(\zb\eta) \d\zb\eta
,\qquad\zb\xi\in\S^{d-1},
\end{equation*}
where $\delta^{(j)}$ denotes the $j$-th derivative of the Dirac delta function and $j\in\N$.
This definition can be imagined by first taking the $j$-th derivative of $f$ in direction of $\zb\xi$ and then computing the integral along the subsphere perpendicular to $\zb\xi$, which is a great circle for $d=3$.
For $j=0$, we obtain the Funk or Funk--Radon transform $\mathcal S^{(0)}_d$, which assigns to $f$ its integrals along all great circles.
The generalized Funk--Radon transform $\mathcal S^{(j)}_d$ belongs to the class of convolution operators on the sphere, cf.\ \cite{HiQu15}.

The aim of the present paper is twofold.
Firstly, we perform a comprehensive analysis on the properties of the generalized Funk--Radon transform $\mathcal S^{(j)}_d$.
Secondly, we utilize Grangeat's formula and our previous findings in order to obtain the singular value decomposition of the cone-beam transform $\mathcal D$.

In \prettyref{thm:Sj-ev-d}, we derive the singular value decomposition of the generalized Funk--Radon transform $\mathcal S^{(j)}_d$. 
This allows us to provide a characterization of its nullspace and range in \prettyref{sec:sobolev}.
In particular, we show that 
$\mathcal S^{(j)}_d$ is a continuous and open operator between the Sobolev spaces $H^s(\S^{d-1}) \to H^{s-j+\frac{d-2}2}(\S^{d-1})$ for $s\in\R$.
This behavior is explained by the fact that $\mathcal S^{(j)}_d$ consists of $j$ derivatives and an integration along a $d-2$ dimensional manifold.
In \prettyref{sec:special-cases-j}, we consider special cases of $\mathcal S^{(j)}$, which equals the hemispherical transform for $j=-1$ and the spherical cosine transform for $j=-2$.
The similarly defined integro-differential Radon transform from Makai et al.\ \cite{MaMaOd00} and the Blaschke--Levy representation coincide with $\mathcal S^{(j)}_d$ for certain but not all parameters, see \prettyref{sec:Sj-similar-transforms}.

Grangeat's formula states a connection between the cone-beam transform $\mathcal D$, the generalized Funk--Radon transform $\mathcal S^{(d-2)}_d$ and the Radon transform, which computes the integrals along hyperplanes in $\R^d$.
Based on Grangeat's formula and our results on $\mathcal S^{(j)}_d$, we derive a singular value decomposition of the cone-beam transform $\mathcal D$ in \prettyref{sec:cone-beam-svd}, where we assume that the function $f$ is supported on the unit ball $\B^d$, the spatial dimension~$d$ is odd and the scanning set is $\Gamma=\S^{d-1}$.
This contains an alternative proof of the result of Kazantsev \cite{Kaz15} for $d=3$.
We analyze the asymptotic behavior of the singular values of the cone-beam transform in \prettyref{sec:cone-beam-sv-bounds}. 
The smallest singular values grow of the order $\mathcal O(m^{-1/2})$ independently of the spatial dimension $d$, which means that the inversion is ill-posed of degree $1/2$, which is the same as for the Radon transform in 2D.
We also obtain a constant upper bound on the singular values.
However, it is open whether this bound is strict.

The outline of this paper is as follows.
In \prettyref{sec:spherical-harmonics}, we summarize some basic facts about spherical harmonics.
In \prettyref{sec:convolution}, we show the singular value decomposition of the generalized Funk--Radon transform $\mathcal S^{(j)}_d$ and its bijectivity in certain Sobolev spaces.
In \prettyref{sec:cone-beam}, we prove the singular value decomposition of the cone-beam transform $\mathcal D$ and compute bounds on the singular values.
Finally, we state our previous results for the practically most relevant case $\B^3$ in \prettyref{sec:R3}.

\section{Harmonic analysis on the sphere}
\label{sec:spherical-harmonics}

In this section, we are going to summarize some basic facts about harmonic analysis
on the $(d-1)$-dimensional unit sphere $\S^{d-1} = \{\zb x\in\R^{d}:\norm{\zb x}=1 \}$ as it can be found in \cite{AtHa12}.
We denote the volume of $\S^{d-1}$ by
\begin{equation}
\abs{\S^{d-1}}
= \int_{\S^{d-1}} \d\S^{d-1}
= \frac{2\pi^{d/2}}{\Gamma(d/2)}
=\begin{cases}
\frac{2\pi^{d/2}}{(\frac{d}{2}-1)!},  & d$ even$\\
\frac{2^{\frac{d+1}{2}} \pi^{\frac{d-1}{2}}}{(d-2)!!}, & d$ odd$.
\end{cases}
\label{eq:VSd}
\end{equation}
In the following, we give an orthonormal basis of the Lebesgue space $L^{2}(\S^{d-1})$, which is the space of square-integrable functions $f\colon\S^{d-1}\to\C$ with the inner product 
$$\left<f,g\right>_{L^2(\S^{d-1})}=\int_{\S^{d-1}} f(\zb\xi)\,g(\zb\xi)\d\zb\xi.$$
The Legendre polynomial $P_{n,d}$ of degree $n\in\N$ in dimension $d\ge3$ is given by the Rodrigues formula \cite[(2.70)]{AtHa12}
\begin{equation}
P_{n,d}(t)
= (-1)^n\, \frac{(d-3)!!}{(2n+d-3)!!}\, (1-t^2)^{\frac{3-d}{2}} \left(\frac{\dx}{\dx t}\right)^n (1-t^2)^{n+\frac{d-3}{2}}
,\qquad t\in[-1,1].
\label{eq:Pnd-Rodrigues}
\end{equation}
The classical Legendre polynomials are $P_n = P_{n,3}$.
Here, we use this name also in the general situation $d>3$ as in \cite{Mue98,AtHa12}.
The Legendre polynomials are orthogonal with respect to the weight function $w_d(t) = (1-t^2)^{\frac{d-3}{2}}$.
They satisfy the orthogonality relation
\begin{equation}
\left<P_{n,d},P_{m,d}\right>_{w_d}
= \int_{-1}^{1} P_{n,d}(t)\, P_{m,d}(t)\, (1-t^2)^{\frac{d-3}{2}} \d t
= \delta_{n,m} \frac{\abs{\S^{d-1}}}{N_{n,d}\abs{\S^{d-2}}}
,
\label{eq:Legendre-orthogonal}
\end{equation}
where 
\begin{equation}
N_{n,d} 
= \frac{(2n+d-2)\,(n+d-3)!}{n!\,(d-2)!}.
\label{eq:Nnd}
\end{equation}
Up to normalization, the Legendre polynomial $P_{n,d}$ is equal to the Gegenbauer or ultraspherical polynomial \cite[(2.145)]{AtHa12}
\begin{equation}
C_n^{(\frac{d-2}{2})} = \binom{n+d-3}{n} P_{n,d}
.
\label{eq:Gegenbauer-Legendre}
\end{equation}
The Gegenbauer polynomial $C_n^{(\alpha)}$ for $\alpha>-1/2$ satisfies the explicit expression \cite[22.3]{AbSt64}
\begin{equation}
C_n^{(\alpha)} (t)
= \frac{1}{\Gamma(\alpha)} \sum_{m=0}^{\left\lfloor \frac{n}{2}\right\rfloor} \frac{(-1)^m\,\Gamma(n-m+\alpha)}{m!\, (n-2m)!} (2t)^{n-2m}.
\label{eq:Gegenbauer-explicit}
\end{equation}
and the Rodrigues formula
\begin{equation}
C_{n}^{(\alpha)}(t)
=\frac{(-1)^n\,\Gamma(\alpha+\frac12)\,\Gamma(n+2\alpha)}{2^n\,n!\,\Gamma(2\alpha)\,\Gamma(\alpha+n+\frac12)}\, (1-t^2)^{\frac12-\alpha} \left(\frac{\dx}{\dx t}\right)^n (1-t^2)^{n+\alpha-\frac{1}{2}}
.
\label{eq:Gegenbauer-Rodrigues}
\end{equation}
We define the space $\mathscr Y_{n,d}(\S^{d-1})$ as the range of the unnormalized projection operator 
$$
L^2(\S^{d-1})\to L^2(\S^{d-1}),\;
f\mapsto \int_{\S^{d-1}} f(\zb\xi)\, P_{n,d}(\zb\xi^\top(\zb\cdot)) \d\zb\xi.
$$
The space $\mathscr Y_{n,d}(\S^{d-1})$ consists of harmonic polynomials that are homogeneous of degree $n$, restricted to the sphere $\S^{d-1}$.
Let 
$Y_{n,d}^k$ for ${k=1,\dots,N_{n,d}},$
be an orthonormal basis of $\mathscr Y_{n,d}(\S^{d-1})$ in $L^2(\S^{d-1})$.
The addition theorem \cite[(2.24)]{AtHa12} states that for $n\in\N$
\begin{equation}
\label{eq:addition-theorem}
\sum_{k=1}^{N_{n,d}} Y_{n,d}^k(\zb\xi)\,\overline {Y_{n,d}^k(\zb\eta)}
= \frac{N_{n,d}}{\abs{\S^{d-1}}}\, P_{n,d}(\zb\xi^\top\zb\eta)
,\qquad \zb\xi,\zb\eta\in\S^{d-1}.
\end{equation}
We define the Gaunt coefficients
\begin{equation*}
G^{n,k,d}_{n_1,k_1,n_2,k_2}
= \int_{\S^{d-1}} Y_{n_1,d}^{k_1}(\zb\xi)\, Y_{n_2,d}^{k_2}(\zb\xi)\, \overline{Y_{n,d}^{k}(\zb\xi)} \d\zb\xi.
\end{equation*}
Then the product of two spherical harmonics can be written as the sum
\begin{equation}
Y_{n_1,d}^{k_1}(\zb\xi)\, Y_{n_2,d}^{k_2}(\zb\xi)
= \sum_{\substack{n=\abs{n_1-n_2}\\n-n_1-n_2\text{ even}}}^{n_1+n_2} \sum_{k=1}^{N_{n,d}} G^{n,k_1+k_2,d}_{n_1,k_1,n_2,k_2}\, Y_{n,d}^{k}
(\zb\xi)
,\qquad \zb\xi\in\S^{d-1}
.
\label{eq:Gaunt-series-d}
\end{equation}

\section{Spherical convolution}
\label{sec:convolution}

The spherical convolution of a function $\psi \colon [-1,1]\to\C$ with a function $f\colon \S^{d-1}\to\C$ is defined by
\begin{equation*}
[\psi\star f](\zb \xi) = \int_{\S^{d-1}} f(\zb{\eta})\, \psi(\zb\xi^\top\zb\eta) \d\zb{\eta}
,\qquad \zb{\xi} \in\S^{d-1}.
\end{equation*}
The Funk--Hecke formula \cite[Thm.\ 2.22]{AtHa12} states that for a spherical harmonic $Y_{n,d}\in\mathscr Y_{n,d}(\S^{d-1})$ and $\psi\in L^1[-1,1]$ where $\int_{-1}^{1} \abs{\psi(t)}\,(1-t^2)^{\frac{d-3}{2}}\d t$ is finite, we have
\begin{equation}
[\psi\star Y_{n,d}](\zb\xi)
= Y_{n,d}(\zb\xi)\,\abs{\S^{d-2}}\,\int_{-1}^{1} \psi(t)\, P_{n,d}(t)\, (1-t^2)^{\frac{d-3}{2}} \d t.
\label{eq:Funk-Hecke}
\end{equation}

\subsection{Generalized Funk--Radon transform}

For $j\in\N$, we define the generalized Funk--Radon transform $\mathcal S^{(j)}_d$  for $f\in C^\infty(\S^{d-1})$ by \cite{Lou16}
\begin{equation*}
\mathcal S^{(j)}_d f (\zb \xi)
= \int_{\S^{d-1}} \delta^{(j)}(\zb\xi^\top\zb\eta)\, f(\zb\eta) \d\zb\eta
,\qquad\zb\xi\in\S^{d-1}.
\end{equation*}
Here, $\delta^{(j)}$ denotes the $j$-th derivative of the Dirac delta distribution, which is defined by its application to a test function $\psi\in C^\infty[-1,1]$
\begin{equation*}
\int_{-1}^{1} \delta^{(j)}(t)\, \psi(t)\d t
= (-1)^j \int_{-1}^{1} \delta(t)\, \psi^{(j)}(t)\d t
= (-1)^j \psi^{(j)}(0).
\end{equation*}
For $j=0$, the operator $\mathcal S^{(0)}_d$ is the Funk--Radon transform, cf.\ \cite[L.\ 2.2]{Rip11}.
As for now, we define $\mathcal S^{(j)}_d$ only for smooth functions.
However, we will later extend it by density to appropriate Sobolev spaces in \prettyref{sec:sobolev}.

\begin{Remark}
	\label{rem:convolution-with-distribution}
	We explain the above definition of $\mathcal S^{(j)}_d$ and justify why we can apply the Funk--Hecke formula \eqref{eq:Funk-Hecke} for $\psi=\delta^{(j)}$.
	Let $f\in C^\infty(\S^{d-1})$.
	We observe that for any $\zb\xi\in\S^{d-1}$
	\begin{equation*}
	\int_{\S^{d-1}} f(\zb\eta) \d\zb\eta
	= \int_{-1}^{1} \frac{1}{\sqrt{1-t^2}} \int_{\zb\xi^\top\zb\eta=t} f(\zb\eta) \d\lambda(\zb\eta) \d t,
	\end{equation*}
	where $\d\lambda$ is the standard surface measure on the sub-sphere $\{\zb\eta\in\S^{d-1}:\zb\xi^\top\zb\eta=t\}$.
	Then we have
	\begin{align*}
	\mathcal S^{(j)}_df (\zb\xi)
	&= \int_{\S^{d-1}} \delta^{(j)}(\zb\xi^\top\zb\eta)\, f(\zb\eta) \d\zb\eta
	\\
	&= \int_{-1}^{1} \delta^{(j)}(t)\, (1-t^2)^{-\frac{1}{2}} \int_{\zb\xi^\top\zb\eta = t} f(\zb\eta)  \d\lambda(\zb\eta)\d t
	\\
	&= (-1)^j \left.\left(\frac{\dx}{\dx t}\right)^j (1-t^2)^{-\frac{1}{2}} \int_{\zb\xi^\top\zb\eta = t} f(\zb\eta)  \d\lambda(\zb\eta) \right|_{t=0}.
	\end{align*}
	We use the following generalized Funk--Hecke formula \cite[(4.2.10)]{BeBuPa68} 
	\begin{equation*}
	\int_{\zb\xi^\top\zb\eta = t} Y_{n,d}^k(\zb\eta) \d \lambda(\zb\eta)
	= \abs{\S^{d-2}}\, (1-t^2)^{\frac{d-2}{2}}\,
	P_{n,d} (t)\, Y_{n,d}^k(\zb\xi)
	,\qquad \zb\xi\in\S^{d-1},\;t\in(-1,1).
	\end{equation*}
	Hence,
	\begin{equation}
	\mathcal S^{(j)}_d Y_{n,d}^k (\zb\xi)
	= \abs{\S^{d-2}}\, (-1)^j\, \left.\left(\frac{\dx}{\dx t}\right)^j (1-t^2)^{\frac{d-3}{2}} P_{n,d}(t) \right|_{t=0}\, Y_{n,d}^k(\zb\xi).
	\label{eq:Sj-Yn}
	\end{equation}
	\eqref{eq:Sj-Yn} can also be obtained by applying the Funk--Hecke formula \eqref{eq:Funk-Hecke} for $\psi=\delta^{(j)}$.
\end{Remark}

In the following, we use double factorials defined by
$n!!=n(n-2)\cdots 2$ for $n$ even or $n!!=n(n-2)\cdots 1$ for $n$ odd and $0!! = 1$.
The Gamma function is defined for $x>0$ by
$
\Gamma(x) = \int_0^\infty y^{x-1}\,\e^{-y}\d y
$
and satisfies $\Gamma(x+1)=x\,\Gamma(x)$ as well as $\Gamma(n)=(n-1)!$ if $n$ is a positive integer.

\begin{Theorem}
	\label{thm:Sj-ev-d}
	Let $j\in\N$.
	The generalized Funk--Radon transform  $\mathcal S^{(j)}_d\colon C(\S^{d-1})\to C(\S^{d-1})$ satisfies the eigenvalue decomposition
	\begin{equation*}
	\mathcal S^{(j)}_d Y_{n,d}^k 
	= \hat{\mathcal S}^{(j)}_d(n)\, Y_{n,d}^k,
	\qquad n\in\N,\;k=1,\dots,N_{n,d},
	\end{equation*} 
	with the eigenvalues for $n+j$ even and ($n\ge j-d+3$ or $d$ even)
	\begin{align}
	\hat{\mathcal S}^{(j)}_d(n)
	&= 
	\abs{\S^{d-2}}\, (-1)^{\frac{n+j}{2}} \frac{(n+j-1)!!\,(d-3)!!}{(n-j+d-3)!!}
	\label{eq:Sj-ev-d}
	\\&
	= 
	\pi^{\frac{d-2}{2}}\, (-1)^{\frac{n+j}{2}}\, 2^{j+1}\, \frac{\Gamma\left(\frac{n+j+1}{2}\right)}{\Gamma\left(\frac{n-j+d-1}{2}\right)}
	\label{eq:Sj-ev-d-gamma}
	\end{align}
	and otherwise
	\begin{equation*}
	\hat{\mathcal S}^{(j)}_d(n)=0.
	\end{equation*}
\end{Theorem}
\begin{proof}
	Let $n\in\N$ and $k\in\{1,\dots,N_{n,d}\}$.
	By \prettyref{eq:Sj-Yn}, we have
	\begin{equation}
	\mathcal S^{(j)}_d Y_{n,d} (\zb\xi)
	= \abs{\S^{d-2}}\, (-1)^j\,Y_{n,d}^k(\zb\xi)\, \left.\left(\frac{\dx}{\dx t}\right)^j P_{n,d}(t)\, (1-t^2)^{\frac{d-3}{2}}\right|_{t=0}.
	\label{eq:SjYn}
	\end{equation}
	By Rodrigues' formula \eqref{eq:Pnd-Rodrigues}, we have
	\begin{equation}
	\left.\left(\frac{\dx}{\dx t}\right)^j P_{n,d}(t)\, (1-t^2)^{\frac{d-3}{2}}\right|_{t=0}
	= (-1)^n\,\frac{(d-3)!!}{(2n+d-3)!!}\left.\left(\frac{\dx}{\dx t}\right)^{n+j} (1-t^2)^{n+\frac{d-3}{2}}\right|_{t=0}.
	\label{eq:dPnd}
	\end{equation}
	We want to apply the generalized binomial theorem, which states for $a,b,z\in\C$
	\begin{equation}
	(a+b)^z 
	=\sum_{k=0}^\infty \binom{z}{k} a^{z-k}\,b^k,
	\label{eq:binomial-theorem}
	\end{equation}
	where the binomial coefficient of $z\in\C$ and $k\in\N$ is defined by
	\begin{equation}
	\binom{z}{k}
	= \frac{z\,(z-1)\cdots(z-k+1)}{k!}.
	\label{eq:binomial}
	\end{equation}
	The binomial theorem implies
	\begin{equation}
	(1-t^2)^{n+\frac{d-3}{2}}
	= \sum_{k=0}^{\infty} \binom{n+\frac{d-3}2}{k} (-1)^k\, t^{2k}.
	\label{eq:df-binomial-2}
	\end{equation}
	Evaluating the $(n+j)$-th derivative of \eqref{eq:df-binomial-2} at $t=0$ and taking into account that 
	$\left.\left(\frac{\dx}{\dx t}\right)^\ell t^{2k}\right|_{t=0}=(2k)!\, \delta_{\ell,2k}$,
	we obtain
	if $n+j$ is even 
	\begin{equation}
	\begin{split}
	\left.\left(\frac{\dx}{\dx t}\right)^{n+j} (1-t^2)^{n+\frac{d-3}{2}}\right|_{t=0}
	&=\sum_{k=0}^{\infty} \binom{n+\frac{d-3}{2}}{k} (-1)^k \left.\left(\frac{\dx}{\dx t}\right)^{n+j} t^{2k}\right|_{t=0}
	\\&= \binom{n+\frac{d-3}{2}}{\frac{n+j}{2}} (-1)^{\frac{n+j}{2}}\, (n+j)!
	\end{split}
	\label{eq:der1}
	\end{equation}
	and zero otherwise.
	By its definition in \eqref{eq:binomial}, the binomial coefficient $\binom{z}{k}$ is zero if and only if both $z$ is a nonnegative integer and $z<k$.
	Hence, the binomial coefficient $\binom{n+\frac{d-3}{2}}{\frac{n+j}{2}}$ from \eqref{eq:der1} is nonzero if and only if $\frac{d-3}{2}$ is not an integer or $n+\frac{d-3}{2} \ge \frac{n+j}{2}$.
	This condition can be simplified to that $d$ is even or $n\ge j-d+3$.
	Then we have
	\begin{equation}
	\begin{split}
	\binom{n+\frac{d-3}{2}}{\frac{n+j}{2}}
	&= \frac{\left(\frac{2n+d-3}{2}\right)\left(\frac{2n+d-3}{2}-1\right)\cdots\left(\frac{2n+d-3}{2}-\frac{n+j}{2}+1\right)}{\left(\frac{n+j}{2}\right)!}
	\\&= \frac{(2n+d-3)!!}{(n-j+d-3)!!\,(n+j)!!}.
	\end{split}
	\label{eq:binomial-df}
	\end{equation}
	Combining \eqref{eq:dPnd}, \eqref{eq:der1} and \eqref{eq:binomial-df}, we obtain
	\begin{align*}
	&\left.\left(\frac{\dx}{\dx t}\right)^j P_{n,d}(t)\, (1-t^2)^{\frac{d-3}{2}}\right|_{t=0}
	\\
	={}& (-1)^n\,\frac{(d-3)!!}{(2n+d-3)!!}
	\frac{(2n+d-3)!!}{(n+j)!!(n-j+d-3)!!}(-1)^{\frac{n+j}{2}} (n+j)!
	\\
	={}& (-1)^{\frac{n-j}{2}}\,\frac{(n+j-1)!!\,(d-3)!!}{(n-j+d-3)!!}
	\end{align*}
	if $n+j$ is even and ($n\ge j-d+3$ or $d$ even), and zero otherwise. Plugging into \eqref{eq:SjYn} shows \eqref{eq:Sj-ev-d}.
	Inserting the volume \eqref{eq:VSd} of $\S^{d-2}$ into \eqref{eq:Sj-ev-d},
	we have for $n+j$ even and ($n\ge j-d+3$ or $d$ even)
	\begin{align*}
	\hat{\mathcal S}^{(j)}_d(n)
	&= 
	\abs{\S^{d-2}}\, (-1)^{\frac{n-j}{2}} \frac{(n+j-1)!!\,(d-3)!!}{(n-j+d-3)!!}
	\\
	&=
	\frac{2\pi^{\frac{d-1}{2}}}{\Gamma\left(\frac{d-1}{2}\right)}\, (-1)^{\frac{n-j}{2}} \frac{(n+j-1)(n+j-3)\cdots 1}{(n-j+d-3)(n-j+d-5)\cdots(d-1)}
	\\&=
	\frac{2\pi^{\frac{d-1}{2}}}{\Gamma\left(\frac{d-1}{2}\right)}\, (-1)^{\frac{n-j}{2}} \frac{2^{\frac{n+j}{2}} \left(\frac{n+j-1}{2}\right) \left(\frac{n+j-1}{2}-1\right)\cdots\left(\frac12\right)}
	{2^{\frac{n-j}{2}}\,\left(\frac{n-j+d-1}{2}-1\right)\left(\frac{n-j+d-1}{2}-2\right)\cdots\left(\frac{d-1}{2}\right)}
	\\&=
	\frac{2\pi^{\frac{d-1}{2}}}{\Gamma\left(\frac{d-1}{2}\right)}\, (-1)^{\frac{n-j}{2}}\, 2^{j}\,\frac{\Gamma\left(\frac{n+j+1}{2}\right)}{\Gamma\left(\frac12\right)}
	\frac{\Gamma\left(\frac{d-1}{2}\right)}{\Gamma\left(\frac{n-j+d-1}{2}\right)}
	\\&=
	\pi^{\frac{d-2}{2}}\, (-1)^{\frac{n-j}{2}}\, 2^{j+1}\, \frac{\Gamma\left(\frac{n+j+1}{2}\right)}{\Gamma\left(\frac{n-j+d-1}{2}\right)},
	\end{align*}
	which shows \eqref{eq:Sj-ev-d-gamma}.
\end{proof}

\prettyref{thm:Sj-ev-d} traces back to \cite{Min04} for $j=0$ on $\S^2$.
The case $j=0$ and $d$ arbitrary was shown in \cite{BeBuPa68}.

\subsection{$\mathcal S^{(j)}$ in Sobolev spaces}

\label{sec:sobolev}

In this section, we extend $\mathcal S^{(j)}$ to a continuous operator between Sobolev spaces.
The spherical Sobolev space $H^s(\S^{d-1})$ of order $s\in\R$ is defined as completion of $C^\infty(\S^{d-1})$ with respect to the Sobolev norm
\begin{equation}
\norm{f}_{H^{s}(\S^{d-1})}^{2}
=\sum_{n=0}^{\infty} \sum_{k=1}^{N_{n,d}} \left(n+\tfrac{d-2}{2}\right)^{2s}\abs{\left<f,Y_{n,d}^k\right>_{L^2(\S^{d-1})}}^{2},
\label{eq:Hs}
\end{equation}
cf.\ \cite[(3.98)]{AtHa12}.
The spherical harmonics $Y_{n,d}^k$ are dense in $H^s(\S^{d-1})$.
The Sobolev spaces are nested: we have $H^s(\S^{d-1}) \hookrightarrow H^t(\S^{d-1})$ whenever $s>t$.
The space $H^0(\S^{d-1})$ can be identified with $L^2(\S^{d-1})$.
If $s$ is a positive integer, $H^s(\S^{d-1})$ can be imagined as the space of functions defined on $\S^{d-1}$ whose (distributional) derivatives up to order $s$ are in $L^2(\S^{d-1})$.

In the following lemma, we derive an asymptotic approximation of the eigenvalues $\hat{\mathcal S}^{(j)}_d(n)$ from \prettyref{thm:Sj-ev-d}. 
We use the notation of asymptotic equivalence $a(n) \simeq b(n)$ for $n\to\infty$ if $\lim_{n\to\infty} {a(n)}/{b(n)}=1$.

\begin{Lemma}
	\label{lem:Sj-ev-d-approx}
	Let $j\in\N$.
	We have for $n\to\infty$ with $n+j$ even and $n\ge j$
	\begin{equation*}
	\abs{\hat{\mathcal S}^{(j)}_d(n)}
	\simeq n^{j-\frac{d-2}2}\, \pi^{\frac{d-1}{2}}\, 2^{\frac{d}{2}},
	\end{equation*}
\end{Lemma}
\begin{proof}
	Let $n+j$ be even and $n\ge j$.
	We apply Stirling's approximation  
	of the Gamma function 
	\begin{equation*}
	\Gamma(x) \simeq \sqrt{2\pi}\, x^{x-\frac12}\,\e^{-x}
	,\qquad x\to\infty
	\end{equation*}
	to the eigenvalues \eqref{eq:Sj-ev-d-gamma} and obtain for $n\to\infty$
	\begin{align*}
	\abs{\hat{\mathcal S}^{(j)}_d(n)}
	&=\pi^{\frac{d-2}{2}}\, 2^{j+1}\, \frac{\Gamma\left(\frac{n+j+1}{2}\right)}{\Gamma\left(\frac{n-j+d-1}{2}\right)}
	\\
	&\simeq\pi^{\frac{d-2}{2}}\, 2^{j+1}\, \frac{\left(\frac{n+j+1}{2}\right)^{\frac{n+j}{2}}\,\e^{-\frac{n+j+1}{2}}}{\left(\frac{n-j+d-1}{2}\right)^{\frac{n-j+d-2}{2}}\, \e^{-\frac{n-j+d-1}{2}} }
	\\
	&= \pi^{\frac{d-2}{2}}\, 2^{\frac{d}{2}}\, \e^{\frac{d-2}{2}-j}\, \frac{\left(n+j+1\right)^{\frac{n+j}{2}}}{\left(n-j+d-1\right)^{\frac{n-j+d-2}{2}} }
	\\
	&= 2^{\frac{d}{2}}\,\pi^{\frac{d-2}{2}}\,\e^{\frac{d-2}{2}-j} \left(1+\frac{2j-d+2}{n-j+d-1}\right)^{\frac{n}{2}} (n+j+1)^{\frac{j}{2}}\, (n-j+d-1)^{\frac{j+2-d}{2}}.
	\end{align*}
	Considering that $\e^{x}=\lim_{n\to\infty} \left(1+\frac{x}{n}\right)^{n}$, we obtain
	\begin{equation*}
	\abs{\hat{\mathcal S}^{(j)}_d(n)}
	\simeq 2^{\frac{d}{2}}\,\pi^{\frac{d-2}{2}}\,\e^{\frac{d-2}{2}-j}\, \e^{\frac{2j-d+2}{2}}\, n^{j+\frac{2-d}{2}}
	= 2^{\frac{d}{2}}\,\pi^{\frac{d-2}{2}}\, n^{j-\frac{d-2}{2}}
	.\qedhere
	\end{equation*}
\end{proof}

The following mapping property of $\mathcal S^{(j)}_d$ between Sobolev spaces was shown for $j=0$ in \cite[§ 4]{Str81}.

\begin{Theorem}
	\label{cor:Sj-Sobolev}
	Let $s\in\R$ and $j\in\N$.
	The generalized Funk--Radon transform $\mathcal S^{(j)}_d$ extends to a  continuous operator 
	\begin{equation}
	\label{eq:Sj-Sobolev}
	\mathcal S^{(j)}_d\colon { H}^s(\S^{d-1}) \to { H}^{s-j+\frac{d-2}2}(\S^{d-1}).
	\end{equation}
	If $j>\frac{d-2}{2}$, then  $\mathcal S^{(j)}_d\colon L^2(\S^{d-1})\to L^2(\S^{d-1})$ is compact.
	The nullspace of $\mathcal S^{(j)}_d$ is the closed linear span 
	\begin{equation*}
	\overline{\operatorname{span}}\left\{ \mathscr Y_{n,d}(\S^{d-1}): \text{$n+j$ odd or ($n\le j-d+1$ and $d$ odd)}  \right\}.
	\end{equation*}
	If $d$ is odd and $j\ge d-1$, the nullspace of $\mathcal S^{(j)}_d$ comprises the sum of all polynomials of degree up to $j-d+1$ and all odd (even) functions whenever $j$ is even (odd).
	Otherwise, the null-space of $\mathcal S^{(j)}_d$ comprises all odd (even) functions whenever $j$ is even (odd).
\end{Theorem}

\begin{proof}
	By the definition \eqref{eq:Hs} of the Sobolev space,
	$\mathcal S^{(j)}_d$ is continuous if and only if the sequence 
	\begin{equation*}
	n\mapsto \abs{\hat{\mathcal{S}}^{(j)}_d(n)}\left(n+\frac{d-2}{2}\right)^{-j+\frac{d-2}{2}}
	\end{equation*} 
	has an upper bound, which follows from \prettyref{lem:Sj-ev-d-approx}.
	The compactness for $j>\frac{d-2}{2}$ follows because then the eigenvalues $\hat{\mathcal S}^{(j)}_d(n)$ converge to 0 for $n\to\infty$.
	The nullspace of $\mathcal S^{(j)}_d$ consists of the closed span of all spherical harmonics $Y_{n,d}^k$ where $n\in\N$ satisfies $\hat{\mathcal{S}}^{(j)}_d(n)=0$.
\end{proof}

The order of smoothness $s-j+\frac{d-2}{2}$ of the Sobolev space in \eqref{eq:Sj-Sobolev} is not unexpected, because $\mathcal S^{(j)}_d$ consists of $j$ differentiations, which lower the order of smoothness by $j$, and the integration along a $(d-2)$-dimensional submanifold, which raises the order of smoothness by $\frac{d-2}{2}$.

\subsection{Special cases of $j$}
\label{sec:special-cases-j}
In this section, we take a look at $\mathcal S^{(j)}_d$ for certain special choices of $j$, some of which are already well-known operators from literature.
Even though $\mathcal S^{(j)}_d$ was initially defined only for $j\in\N$, we can extend it to negative $j$ by the singular value decomposition \eqref{eq:Sj-ev-d}.
Inserting $j=-1$ in equation \eqref{eq:Sj-ev-d-gamma} of the eigenvalues yields for odd $n$
\begin{equation*}
\hat{\mathcal S}^{(-1)}_d(n)
= \abs{\S^{d-2}} (-1)^{\frac{n-1}{2}} \frac{(n-2)!!\,(d-3)!!}{(n+d-2)!!}
= 2\,(-1)^{\frac{n-1}{2}}\, \frac{2\pi^{\frac{d-1}{2}}}{\Gamma\left(\frac{1}{2}\right)}\, \frac{\Gamma\left(\frac{n}{2}\right)}{\Gamma\left(\frac{n+d}{2}\right)}.
\end{equation*}
Hence, ${\mathcal S}^{(-1)}_d$ is the modified hemispherical transform \cite{Rub99}
\begin{equation}
\mathcal S^{(-1)}_d f(\zb\xi) = \frac12 \int_{\S^{d-1}} \operatorname{sgn(\zb\xi^\top\zb\eta)}\,f(\zb\eta)\d\zb\eta.
\label{eq:hemispherical}
\end{equation}
Inserting $j=-2$ gives the eigenvalues
$$
\hat{\mathcal S}^{(-2)}_d (n) 
= \begin{cases}
2 \abs{\S^{d-2}} (-1)^{\frac{n-2}{2}} \frac{(n-3)!!\,(d-2)!!}{(n+d-1)!!}
,&\text{$n$ even}\\
0,&\text{$n$ odd}
\end{cases}
$$
of the spherical cosine transform, cf.\ \cite[L.\ 3.4.5]{Gro96}, 
$$
\mathcal S^{(-2)}_d f(\zb\xi) = \frac12 \int_{\S^{d-1}} \abs{\zb\xi^\top\zb\eta}\,f(\zb\eta)\d\zb\eta.
$$
In the case $d = 2j+2$, which in particular implies that $d$ is even, we have the eigenvalues
\begin{equation*}
\hat{\mathcal S}^{(j)}_{2j+2} (n) 
= \begin{cases}
\abs{\S^{2j}}\, (-1)^{\frac{n+j}{2}}\,(2j-1)!!=
2\,(2\pi)^j\,(-1)^{\frac{n+j}{2}}
,&\text{$n+j$ even}\\
0,&\text{$n+j$ odd},
\end{cases}
\end{equation*}
which are, except for their sign, independent of $n$.
Hence, if $j$ is even (odd), the operator ${\mathcal S}^{(j)}_{2j+2}\colon L^2(\S^{2j+1})\to L^2(\S^{2j+1})$ restricted to the even (odd) functions is an isometry.

The following theorem shows an inversion formula for the Funk--Radon transform in even dimensions.

\begin{Theorem}
	Let $d\ge2$ be even.
	Then any even function $f\colon\S^{d-1}\to\C$ can be reconstructed from its Funk--Radon transform $g = \mathcal S^{(0)}_d f$ by
	$$
	f = \frac{1}{\abs{\S^{d-2}}^2\,((d-3)!!)^2}\, \mathcal S^{(d-2)}_d g.
	$$
\end{Theorem}
\begin{proof}
	Let $n\in\N$ be even.
	On the one hand, we have the eigenvalues
	\begin{equation*}
	\hat{\mathcal S}^{(d-2)}_{d} (n) = 
	\abs{\S^{d-2}}\, (-1)^{\frac{n+d-2}{2}} \frac{(n+d-3)!!\,(d-3)!!}{(n-1)!!}.
	\end{equation*}
	On the other hand, the Funk--Radon transform $\mathcal S^{(0)}_d$ has the eigenvalues
	\begin{equation*}
	\hat{\mathcal S}^{(0)}_{d} (n)  =
	\abs{\S^{d-2}}\, (-1)^{\frac{n}{2}} \frac{(n-1)!!\,(d-3)!!}{(n+d-3)!!}.
	\end{equation*}
	Hence, the product of the two operators has the constant eigenvalues
	\begin{equation*}
	\widehat{\mathcal S^{(d-2)}_{d}\mathcal S^{(0)}_{d} } (n) =
	\abs{\S^{d-2}}^2\, (-1)^{\frac{d-2}{2}}\,(d-3)!!^2.\qedhere
	\end{equation*}
\end{proof}

\subsection{Similar transforms}
\label{sec:Sj-similar-transforms}
In this section, we consider two integral transforms, which are equal to $\mathcal S^{(j)}_d$ for certain but not all parameters $j$.

\paragraph{Integro-differential transform}
For $j\in \N$ and $\vartheta\in\left[-\frac\pi2,\frac\pi2\right]$, we define the integro-differential transform $\mathcal R^{(j)}_{d,\vartheta}\colon C^j(\S^{d-1})\to C(\S^{d-1})$ by
\begin{equation*}
\mathcal R^{(j)}_{d,\vartheta} f(\zb\xi) = \int_{\zb\xi^\top\zb\omega = 0} 
\left(\frac{\partial}{\partial \vartheta}\right)^j f(\zb\xi\sin\vartheta + \zb\omega \cos\vartheta) \d\zb\omega,
\qquad\zb\xi\in \S^d,
\end{equation*}
which was introduced in \cite{MaMaOd00}.
The operator $\mathcal R^{(0)}_{d,\vartheta}$ has been investigated in \cite{Sch69}.
For $j>0$, we first take the $j$-th derivative of $f$ perpendicular to the circle of integration.
It was shown in \cite{MaMaOd01} that the nullspace of the operator $\mathcal R^{(j)}_{d,\vartheta}$ is 
\begin{equation*}
\overline{\operatorname{span}}\left\{\mathscr Y_{n,d}(\S^{d-1}): \left(\frac{\dx}{\dx \vartheta}\right)^j P_{n,d}(\sin\vartheta)=0 \right\}.
\end{equation*}
If $\vartheta=0$, we write $\mathcal R^{(j)}_d = \mathcal R^{(j)}_{d,0}$.
If $j\ge1$, the null–space of $\mathcal R^{(j)}_d \colon C^j(\S^{d-1}) \to C(\S^{d-1})$ equals for $j$ odd (even) the set $\{f \in C^m(\S^d) : f$ is even ($f$ is the sum of an odd function and a constant)$\}$.

\begin{Theorem}
	\label{thm:SR}
	We have $\mathcal S^{(0)}_d = \mathcal R^{(0)}_d$ and $\mathcal S^{(1)}_d = -\mathcal R^{(1)}_d$.
\end{Theorem}
\begin{proof} 
	For $j=0$, we see that $\mathcal S^{(0)}_d = \mathcal R^{(0)}_d$ is the Funk--Radon transform.
	By \cite{MaMaOd01}, the operator $\mathcal R^{(j)}_d$ has as eigenfunctions the spherical harmonics $Y_{n,d}^k$ and the eigenvalues 
	$$\hat{\mathcal R}^{(j)}_d(n) 
	= \abs{\S^{d-2}} \left.\frac{\mathrm d}{\mathrm d\vartheta} P_{n,d}(\sin\vartheta)\right|_{\vartheta=0}.$$
	\prettyref{thm:Sj-ev-d} shows that $\mathcal S^{(j)}_d$ has the same eigenfunctions.
	Hence, the two operators coincide if their respective eigenvalues do.
	We have for $j=1$ on the one hand
	\begin{equation*}
	\hat{\mathcal R}^{(1)}_d (n)
	= \abs{\S^{d-2}} \left.P_{n,d}'(\sin\vartheta)\,\cos\vartheta\right|_{\vartheta=0}
	= \abs{\S^{d-2}} P_{n,d}'(0)
	\end{equation*}
	and on the other hand
	\begin{equation*}
	\hat{\mathcal S}^{(1)}_d(n)
	= -\abs{\S^{d-2}}\left.\frac{\dx}{\dx t} P_{n,d}(t)\, (1-t^2)^{\frac{d-3}{2}}\right|_{t=0}
	= -\abs{\S^{d-2}}P_{n,d}'(0).\qedhere
	\end{equation*}
\end{proof}
\begin{Remark}
	\prettyref{thm:SR} does not hold for all $j$.
	We have for $d=3$ and $j=2$
	\begin{equation*}
	\hat{\mathcal R}^{(2)}_3(n) 
	= 2\pi \big(P_n''(\sin\vartheta)\,\cos^2\vartheta-P_n'(\sin\vartheta)\,\sin\vartheta\big)\Big|_{\vartheta=0}
	= 2\pi\,P_n''(0) = \hat{\mathcal S}^{(2)}_{3}(n).
	\end{equation*}
	However, for $j=3$
	\begin{align*}
	\hat{\mathcal R}^{(3)}_3(n) 
	&= 2\pi \big(P_n'''(\sin\vartheta)\,\cos^3\vartheta -3P_n''(\sin\vartheta)\,\cos\vartheta\,\sin\vartheta -P_n'(\sin\vartheta)\,\cos\vartheta\big)\Big|_{\vartheta=0}
	\\&= 2\pi\,(P_n'''(0)-P_n'(0))
	\end{align*}
	does not coincide with $-\hat{\mathcal S}^{(3)}_3(n) = 2\pi\,P_n'''(0)$.
\end{Remark}

\paragraph{Blaschke--Levy reprsentation}
Another related transform is the $\alpha$-cosine transform or Blaschke--Levy representation \cite{Kol97,Rub99}
\begin{equation*}
\mathcal H^{(\alpha)}f(\zb\xi)
= \int_{\S^{d-1}} \abs{\zb\xi^\top\zb\eta}^\alpha f(\zb\eta)\d\zb\eta
\end{equation*}
with singular values 
\begin{equation*}
{\hat {\mathcal H}^{(\alpha)}}(n)
=
\begin{cases} \displaystyle(-1)^{n/2}\frac{\Gamma\left(\frac{n-\alpha}{2}\right)}{\Gamma\left(\frac{n+d+\alpha}{2}\right)}, &n\text{ even}\\
\displaystyle 0, &n\text{ odd}.
\end{cases}
\end{equation*}
Hence, for $j$ even and $\alpha=-j-1$, the $\alpha$-cosine transform $\mathcal H^{\alpha}$ is, up to a constant factor, equal to the generalized Funk--Radon transform $\mathcal S^{(j)}_{d}$.

\section{Cone-beam transform}

\label{sec:cone-beam}

\subsection{Connection of Radon and cone-beam transform}

\paragraph{Radon transform}
We define the Radon transform $\mathcal R$ on the $d$-dimensional unit ball $\B^d = \{\zb x\in\R^d: \norm{x}\le1 \}$ by
\cite[Sec.~II.1]{Natterer86}
\begin{equation}
\begin{split}
&\mathcal R\colon L^2(\B^d) \to L^2(\S^{d-1}\times[-1,1],w_{d/2}^{-1})
\\&
\mathcal Rf(\zb\omega, s)
= \int_{\zb x^\top \zb\omega=s} f(\zb x) \d\zb x
\end{split}
\end{equation}
with the weight function 
$$ w_\nu(s) = (1-s^2)^{\nu-1/2}
,\qquad s\in[-1,1].$$
The Radon transform on the unit ball has the following singular value decomposition~\cite{Lou84}.
For $m\in\N$, $l=0,\dots,m$ with $m+l$ even and $k=1,\dots,N_{l,d}$, we have
\begin{equation}
\mathcal R \widetilde V_{m,l,k}(\zb\omega,s)
= \frac{\sqrt{2m+d}\,\Gamma(\frac{d}{2})\, m!}{2^{1-d}\, \pi^{1-\frac{d}{2}}\, (m+d-1)!}\, (1-s^2)^{\frac{d-1}{2}}\, C_m^{(\frac{d}{2})}(s)\, Y_{l,d}^k(\zb\omega),
\label{eq:Radon-svd}
\end{equation}
where
\begin{equation}
\widetilde V_{m,l,k}(s\zb\omega)
= \sqrt{2m+d}\, s^l P_{\frac{m-l}{2}}^{\left(0, l+\frac{d-2}{2}\right)} (2s^2-1)\, Y_{l,d}^k(\zb\omega)
,\qquad s\in[0,1],\; \zb\omega\in\S^{d-1}
\label{eq:Vmlk}
\end{equation}
and $P_n^{(\alpha,\beta)}$ denotes the Jacobi polynomial of degree $n$ and orders $\alpha,\beta>-1$.
The set 
$$\left\{\widetilde V_{m,l,k}: l\in\N,\; m\in\{l,l+2,l+4,\dots\},\; k\in\{1,\dots,N_{l,d}\}\right\}$$ 
is an orthonormal basis of $L^2(\B^d)$ consisting of polynomials of degree $m\in\N$.

\begin{Remark}
	We compute the norm of $\widetilde V_{m,l,k}$ in $L^2(\B^d)$.
	The Jacobi polynomials satisfy
	\begin{equation*}
	\int_{-1}^1 (1-x)^{\alpha} (1+x)^{\beta} P_m^{(\alpha,\beta)} (x)P_n^{(\alpha,\beta)} (x)\,dx =\frac{2^{\alpha+\beta+1}}{2n+\alpha+\beta+1} \frac{\Gamma(n+\alpha+1)\Gamma(n+\beta+1)}{\Gamma(n+\alpha+\beta+1)n!}.
	\end{equation*}
	We have
	\begin{equation*}
	\norm{\widetilde V_{m,l,k}}_{L^2(\B^d)}^2
	= (2m+d)
	\int_{0}^{1} s^{2l} \abs{P_{\frac{m-l}{2}}^{\left(0, l+\frac{d-2}{2}\right)} (2s^2-1)}^2 s^{d-1} \d s 
	\int_{\S^{d-1}} \abs{Y_{l,d}^k(\zb\omega)}^2 \d\zb\omega.
	\end{equation*}
	By the normalization of the spherical harmonics and the substitution $t=2s^2-1$ with $\d t = 4s\d s$, we obtain
	\begin{align*}
	\norm{\widetilde V_{m,l,k}}_{L^2(\B^d)}^2
	&= \frac{2m+d}{4}
	\int_{-1}^{1} \left(\frac{t+1}{2}\right)^{l+\frac{d-2}{2}} \abs{P_{\frac{m-l}{2}}^{\left(0, l+\frac{d-2}{2}\right)} (t)}^2 \d t
	\\
	&= \frac{2m+d}{4} 2^{-l-\frac{d-2}{2}} \frac{2^{l+\frac{d}{2}}}{m+\frac{d}{2}} \frac{\Gamma(\frac{m-l}{2}+1)\, \Gamma(\frac{m+l+d}{2})}{(\frac{m-l}{2})!\, \Gamma(\frac{m+l+d}{2})}
	= 1
	.
	\end{align*}
\end{Remark}

\paragraph{Cone-beam transform}

The cone-beam transform, which is also known as divergent beam X-ray transform, with scanning set $\Gamma\subset\R^d$ is defined by
\begin{equation*}
\mathcal Df (\zb a,\zb\omega)
= \int_{0}^{\infty} f(\zb a + t\zb\omega)\d t
,\qquad \zb\omega\in\S^{d-1},\; \zb a\in\Gamma.
\end{equation*}
\paragraph{Grangeat's formula}
There is a relation between the Radon transform and the cone-beam transform.
Let $h\colon\R\to\R$ be a function that is homogeneous of degree $1-d$. It was essentially shown in \cite{HaSmSoWa80} (see also \cite[Sec.\ 2.3]{NaWue00} and \cite[Sec.~2.2.1]{Pal16}) that
\begin{equation}
\int_{-\infty}^{\infty} \mathcal Rf(\zb\omega, s)\, h(s-\zb a^\top\zb\omega) \d s
= \int_{\S^{d-1}} \mathcal Df(\zb a,\zb\xi)\, h(\zb\omega^\top\zb\xi) \d\zb\xi.
\label{eq:Hamaker}
\end{equation}
Inserting $h=\delta^{(d-2)}$, we obtain
Grangeat's formula, which was originally proved for $d=3$ in \cite{Gra91}, stating that
\begin{equation}
(-1)^{d}\left(\frac{\partial}{\partial s}\right)^{d-2} \mathcal Rf(\zb\omega, \zb a^\top\zb\omega)
= \mathcal S^{(d-2)}_d \mathcal Df(\zb a,\zb\omega),
\label{eq:Grangeat}
\end{equation}
where $\mathcal S^{(d-2)}_d$ is applied with respect to $\zb\omega$.

The following theorem gives an alternative version of Grangeat's formula for $d=3$.
However, it is not a special case of \eqref{eq:Hamaker}, because the function $h$ is not homogeneous of degree $-2$.

\begin{Theorem}
	\label{thm:Grangeat-variant}
	Let $\zb\omega\in\S^2$ and $\zb a\in\R^3$. 
	We have
	\begin{equation}
	- \mathcal S^{(-1)}_3 \frac{\partial}{\partial s} \mathcal Rf(\zb\omega, \zb a^\top\zb\omega)
	= \int_{\S^2} h(\zb\xi^\top\zb\omega)\, \mathcal Df(\zb a,\zb\xi) \d\zb\xi,
	\label{eq:Grangeat-variant}
	\end{equation}
	where 	
	\begin{equation*}
	h(x) = \frac{2x}{\sqrt{1-x^2}}
	,\qquad x\in(-1,1).
	\end{equation*}
	and $\mathcal S^{(-1)}_{d}$ is the modified hemispherical transform \eqref{eq:hemispherical} applied with respect to $\zb\omega$.
\end{Theorem}
\begin{proof}
	Multiplying Grangeat's formula \eqref{eq:Grangeat} for $d=3$ with $\mathcal S^{(-1)}$, we obtain
	\begin{equation}
	- \mathcal S^{(-1)}_3 \frac{\partial}{\partial s} \mathcal Rf(\zb\omega, \zb a^\top\zb\omega)
	= \mathcal S^{(-1)}_{3} \mathcal S^{(1)}_{3} \mathcal Df(\zb a,\zb\omega).
	\label{eq:Grangeat-variant1}
	\end{equation}
	The left-hand side of the previous equation \eqref{eq:Grangeat-variant1} is the same as that of equation \eqref{eq:Grangeat-variant}.
	We are going to show the equality of the right-hand side of \eqref{eq:Grangeat-variant} and \eqref{eq:Grangeat-variant1} by evaluation for all spherical harmonics $\mathcal Df = Y_n^k$, the assertion for general function $\mathcal Df$ follows by the density of the spherical harmonics $Y_n^k$.
	By \cite[8.922.4]{GrRy07}, we have
	\begin{equation*}
	h(x) 
	= \frac{2x}{\sqrt{1-x^2}}
	= \pi\, \sum_{\substack{n=1\\n\text{ odd}}}^{\infty} {(2n+1)} \frac{(n-2)!!\,n!!}{(n-1)!!\, (n+1)!!}\, P_n(x)
	,\qquad x\in(-1,1).
	\end{equation*}
	Then we have by the Funk--Hecke formula \eqref{eq:Funk-Hecke} for all $n\in\N$, $\abs{k}\le n$ and $\zb\xi\in\S^2$
	\begin{align*}
	\int_{\S^2} h(\zb\xi^\top\zb\omega)\, Y_n^k(\zb\xi) \d\zb\xi
	&= 2\pi \int_{-1}^1 h(t)\, P_n(t) \d t \, Y_n^k(\zb\omega)
	\\&
	= \begin{cases}4\pi^2 \frac{(n-2)!!\,n!!}{(n-1)!!\, (n+1)!!}\, Y_n^k(\zb\omega), &n\text{ odd}\\
	0, &n\text{ even}.
	\end{cases}
	\end{align*}
	On the other hand, the right-hand side of \eqref{eq:Grangeat-variant1} evaluates for odd $n$ by \eqref{eq:Sj-ev-d}
	\begin{equation*}
	{\mathcal S}^{(-1)}_3\, {\mathcal S}^{(1)}_3\, Y_n^k
	=\hat{\mathcal S}^{(-1)}_3(n)\, \hat{\mathcal S}^{(1)}_3(n)\, Y_n^k
	= 4\pi^2 \frac{(n-2)!!\,n!!}{(n-1)!!\, (n+1)!!}\, Y_n^k,
	\end{equation*}
	which shows the assertion.
\end{proof}

\subsection{Singular value decomposition}
\label{sec:cone-beam-svd}
In the following, we consider the cone-beam transform $\mathcal D$ with scanning set $\Gamma=\S^{d-1}$ and we assume that the function $f$ is supported in the unit ball $\B^d$.
We see that $\mathcal Df(\zb a,\zb\omega)=0$ for all $\zb\omega\in\S^{d-1}$ with $\zb a^\top\zb\omega\ge0$ since the ray of integration is outside $\B^d$.
We denote the odd part of the cone-beam transform $\mathcal D(\zb a,\zb\cdot)$ by
$$\mathcal D^\text{(odd)}f(\zb a,\zb\omega) = \frac{Df(\zb a,\zb\omega) - Df(\zb a,-\zb\omega)}{2}.$$
Then 
$$\mathcal Df(\zb a,\zb\omega) = 2\mathcal D^\text{(odd)}f(\zb a,\zb\omega)$$ 
for all $\zb\omega\in\S^2$ with $\zb a^\top\zb\omega<0$ and $\mathcal Df(\zb a,\zb\omega) =0$ otherwise.

\begin{Lemma}
	Let $m\in\N$ and $d\ge3$ be odd.
	Then
	\begin{equation}
	C_{m+1}^{(\frac{d-2}{2})} (s)
	= (-1)^{\frac{d-1}{2}} \frac{(d-2)\,m!}{(m+d-1)!}\, \left(\frac{\partial}{\partial s}\right)^{d-2} (1-s^2)^{\frac{d-1}{2}}\, C_m^{(\frac{d}{2})}(s).
	\label{eq:Gegenbauer-diff}
	\end{equation}
	\label{lem:Gegenbauer-diff}
\end{Lemma}
\begin{proof}
	We denote the right-hand side of \eqref{eq:Gegenbauer-diff} by $B_{m+1}^{(\frac{d-2}{2})}(s)$, 
	i.e., we set
	\begin{equation*}
	B_{m+1}^{(\frac{d-2}{2})} (s)
	= (-1)^{\frac{d-1}{2}} \frac{(d-2)\,m!}{(m+d-1)!}\, \left(\frac{\partial}{\partial s}\right)^{d-2} (1-s^2)^{\frac{d-1}{2}}\, C_m^{(\frac{d}{2})}(s)
	.
	\end{equation*}
	We obtain with Rodrigues' formula \eqref{eq:Gegenbauer-Rodrigues} for the Gegenbauer polynomials
	\begin{align*}
	B_{m+1}^{(\frac{d-2}{2})} (s)
	&= (-1)^{\frac{d-1}{2}} \frac{(d-2)\,m!}{(m+d-1)!}\, \frac{(-1)^m\, (\frac{d-1}{2})!\,(m+d-1)!}{2^m\,m!\,(d-1)!\,(m+\frac{d-1}{2})!}\, \left(\frac{\partial}{\partial s}\right)^{m+d-2}  (1-s^2)^{m+\frac{d-1}{2}}
	\\
	&= \frac{(-1)^{m+\frac{d-1}{2}}\, (d-2)\, (\frac{d-1}{2})!}{2^m\, (d-1)!\, (m+\frac{d-1}{2})!}\, \left(\frac{\partial}{\partial s}\right)^{m+d-2}  (1-s^2)^{m+\frac{d-1}{2}}.
	\end{align*}
	We compute with the binomial theorem \eqref{eq:binomial-theorem}
	\begin{align*}
	B_{m+1}^{(\frac{d-2}{2})} (s)
	&= \frac{(-1)^{m+\frac{d-1}{2}}\, (d-2)\, (\frac{d-1}{2})!}{2^m\, (d-1)!\, (m+\frac{d-1}{2})!}\, \sum_{i=0}^{m+\frac{d-1}{2}} (-1)^i \binom{m+\frac{d-1}{2}}{i} \left(\frac{\partial}{\partial s}\right)^{m+d-2}  s^{2i}
	\end{align*}
	Considering 
	$$
	\left(\frac{\partial}{\partial s}\right)^{k}  s^{2i} = \frac{(2i)!}{(2i-k)!} s^{2i-k},
	$$
	we have
	\begin{multline*}
	B_{m+1}^{(\frac{d-2}{2})} (s)\\
	= \frac{(-1)^{m+\frac{d-1}{2}}\, (d-2)\, (\frac{d-1}{2})!}{2^m\, (d-1)!\, (m+\frac{d-1}{2})!}\,\sum_{i=\left\lceil \frac{m+d-2}{2} \right\rceil}^{m+\frac{d-1}{2}}  \binom{m+\frac{d-1}{2}}{i} \frac{(-1)^i\, (2i)!}{(2i-m-d+2)!}\, s^{2i-m-d+2}
	.
	\end{multline*}
	Shifting the index $i\mapsto l$ with $i = m-l+\frac{d-1}{2}$, we obtain
	\begin{align*}
	&B_{m+1}^{(\frac{d-2}{2})} (s)
	\\
	&= \frac{(-1)^{m+\frac{d-1}{2}}\, (d-2)\, (\frac{d-1}{2})!}{2^m\, (d-1)!\, (m+\frac{d-1}{2})!}\,
	\sum_{l=0}^{\left\lfloor \frac{m+1}{2}\right\rfloor}  \frac{(-1)^{m-l+\frac{d-1}{2}}\,(m+\frac{d-1}{2})!}{(m-l+\frac{d-1}{2})!\,l!} \frac{(2m-2l+d-1)!}{(m+1-2l)!}\, s^{m+1-2l}
	\\
	&= \frac{(d-2)\, (\frac{d-1}{2})!}{2^m\,(d-1)!}\, 
	\sum_{l=0}^{\left\lfloor \frac{m+1}{2}\right\rfloor}  \frac{(-1)^{l}}{(m-l+\frac{d-1}{2})!\,l!} \frac{(2m-2l+d-1)!}{(m+1-2l)!}\, s^{m+1-2l}.
	\end{align*}
	Because
	\begin{equation*}
	\frac{(2m)!}{m!}
	= \frac{2^m\,(2m)!}{(2m)!!}
	= 2^m\,(2n-1)!!,
	\end{equation*}
	we have
	\begin{align*}
	B_{m+1}^{(\frac{d-2}{2})} (s)
	&= \frac{(d-2)}{2^m\, (d-2)!!\,2^{\frac{d-1}{2}}}\, 
	\sum_{l=0}^{\left\lfloor \frac{m+1}{2}\right\rfloor} \frac{(-1)^{l}\,2^{m-l+\frac{d-1}{2}}\, (2m-2l+d-2)!!}{l!\,(m+1-2l)!}\, s^{m+1-2l}
	\\
	&= \frac{1}{(d-4)!!}\, 
	\sum_{l=0}^{\left\lfloor \frac{m+1}{2}\right\rfloor} \frac{(-1)^{l}\,2^{-l}\, (2m-2l+d-2)!!}{l!\,(m+1-2l)!}\, s^{m+1-2l}.
	\end{align*}
	We rewrite the quotient of double factorials with the Gamma function
	\begin{equation*}
	\frac{(m+2k)!!}{m!!}
	=2^k \left(\frac{m+2k}{2}\right)\left(\frac{m+2k-2}{2}\right)\cdots\left(\frac{m+2}{2}\right)
	=2^k \frac{\Gamma\left(\frac{m+2k+2}{2}\right)}{\Gamma\left(\frac{2k+2}{2}\right)}
	\end{equation*}
	and obtain
	\begin{equation*}
	B_{m+1}^{(\frac{d-2}{2})} (s)
	= \sum_{l=0}^{\left\lfloor \frac{m+1}{2}\right\rfloor} \frac{(-1)^{l}\, \Gamma(m-l+\frac{d}{2})}{\Gamma(\frac{d-2}{2})\,l!\,(m+1-2l)!}\, (2s)^{m+1-2l},
	\end{equation*}
	which is exactly the formula \eqref{eq:Gegenbauer-explicit} for the Gegenbauer polynomial $C_{m+1}^{(\frac{d-2}{2})} (s)$.
\end{proof}

\begin{Theorem}
	\label{thm:cone-svd}
	Let $m\in\N$, $l=0,\dots,m$ with $l+m$ even, $k\in\{1,\dots,N_{l,d}\}$ and $d\ge3$ odd.
	The odd cone-beam transform $\mathcal D^{(odd)}\colon \B^d\to\S^{d-1}\times\S^{d-1}$ satisfies for $\zb a,\zb\omega\in\S^{d-1}$
	\begin{equation*}
	\mathcal D^\text{(odd)} \widetilde V_{m,l,k}(\zb a,\zb\omega)
	= \mu_{m,d}\, \sum_{j=1}^{N_{m+1,d}} \overline{Y_{m+1,d}^j(\zb a)}
	\sideset{}{'}\sum_{n=m+1-l}^{l+m+1}  \nu_{n,d} \sum_{i=1}^{N_{n,d}} G^{n,i,d}_{m+1,j,l,k} Y_{n,d}^{i}(\zb\omega),
	\end{equation*}
	where $\sum'$ denotes the summation over odd indices,
	$\widetilde V_{m,l,k}$ is given in \eqref{eq:Vmlk}
	and
	\begin{align}
	\mu_{m,d}
	&= \sqrt{\frac{2^{d+1}\,\pi^{d-1}}{2m+d}},
	\label{eq:mu}
	\\
	\nu_{n,d}
	&= \frac{(-1)^{\frac{n+1}{2}}\,(n-1)!!}{(n+d-3)!!}.
	\label{eq:nu}
	\end{align}
\end{Theorem}
\begin{proof}
	Let $m\in\N$, $l\in\{0,\dots,m\}$ with $m+l$ even, $k\in\{1,\dots,N_{l,d}\}$ and $d$ be odd.
	We have by the singular value decomposition \eqref{eq:Radon-svd} of the Radon transform and \prettyref{lem:Gegenbauer-diff}
	\begin{align*}
	&\left(\frac{\partial}{\partial s}\right)^{d-2} \mathcal R\widetilde V_{m,l,k}(\zb\omega, s)
	\\
	={}& \frac{2^{d-1}\,\pi^{\frac{d}{2}-1}\,\sqrt{2m+d}\,\Gamma(\frac{d}{2})\, m!}{(m+d-1)!}\,\left(\frac{\partial}{\partial s}\right)^{d-2} (1-s^2)^{\frac{d-1}{2}}\, C_m^{(\frac{d}{2})}(s)\, Y_{l,d}^k(\zb\omega)
	\\
	={}& \frac{2^{d-1}\,\pi^{\frac{d}{2}-1}\,\sqrt{2m+d}\,\Gamma(\frac{d}{2})\, m!}{(m+d-1)!}\,(-1)^{\frac{d-1}{2}}\, \frac{(m+d-1)!}{(d-2)\,m!}\, C_{m+1}^{(\frac{d-2}{2})}(s)\, Y_{l,d}^k(\zb\omega)
	\\
	={}& \frac{2^{d-1}\,\pi^{\frac{d}{2}-1}\,\sqrt{2m+d}\,\Gamma(\frac{d}{2})}{d-2}\, (-1)^{\frac{d-1}{2}}\, C_{m+1}^{(\frac{d-2}{2})}(s)\, Y_{l,d}^k(\zb\omega).
	\end{align*}
	By Grangeat's formula \eqref{eq:Grangeat} and the relation \eqref{eq:Gegenbauer-Legendre} between the Gegenbauer and the Legendre polynomials, we obtain
	\begin{multline*}
	\mathcal S_d^{(d-2)} \mathcal D \widetilde V_{m,l,k}(\zb a,\zb\omega)
	\\
	= (-1)^{\frac{d+1}{2}}
	\frac{2^{d-1}\,\pi^{\frac{d}{2}-1}\,\sqrt{2m+d}\,\Gamma(\frac{d}{2})\,(m+d-2)!}{(m+1)!\,(d-2)!} 
	P_{m+1,d}(\zb a^\top\zb\omega)\, Y_{l,d}^k(\zb\omega).
	\end{multline*}
	By the addition formula \eqref{eq:addition-theorem} for spherical harmonics, we have
	\begin{equation*}
	\mathcal S_d^{(d-2)} \mathcal D\widetilde V_{m,l,k}(\zb a,\zb\omega)
	= (-1)^{\frac{d+1}{2}} \abs{\S^{d-1}} 
	\frac{2^{d-1}\,\pi^{\frac{d}{2}-1}\,\Gamma(\frac{d}{2})}{\sqrt{2m+d}}
	\sum_{j=1}^{N_{m+1,d}} \overline{Y_{m+1,d}^j(\zb a)}\, Y_{m+1,d}^j(\zb\omega)\,  Y_{l,d}^k(\zb\omega).
	\end{equation*}
	By the multiplication formula \eqref{eq:Gaunt-series-d} for spherical harmonics, we see that
	\begin{multline*}
	\mathcal S_d^{(d-2)} \mathcal D\widetilde V_{m,l,k}(\zb a,\zb\omega)
	\\
	= (-1)^{\frac{d-1}{2}} \abs{\S^{d-1}} \frac{2^{d+1}\,\pi^{\frac{d}{2}-1}\,\Gamma(\frac{d}{2})}{\sqrt{2m+d}} \sum_{j=1}^{N_{m+1,d}} \overline{Y_{m+1,d}^j(\zb a)}
	\sideset{}{'}\sum_{n=m+1-l}^{m+1+l} \sum_{i=1}^{N_{n,d}} G^{n,i,d}_{m+1,j,l,k} Y_{n,d}^{i}(\zb\omega).
	\end{multline*}
	Since the generalized Funk--Radon transform $\mathcal S^{(d-2)}$ acts only on odd functions, we have $\mathcal S^{(d-2)}_d \mathcal D = \mathcal S^{(d-2)}_d \mathcal D^\text{(odd)}$.
	Then the eigenvalue decomposition \eqref{eq:Sj-ev-d} of  $\mathcal S_d^{(d-2)}$ yields
	\begin{align*}
	\mathcal D^\text{(odd)}\widetilde V_{m,l,k}(\zb a,\zb\omega)
	={}& (-1)^{\frac{d+1}{2}} \frac{\abs{\S^{d-1}}}{\abs{\S^{d-2}}} \frac{2^{d-1}\,\pi^{\frac{d}{2}-1}\,\Gamma(\frac{d}{2})}{\sqrt{2m+d}} \sum_{j=1}^{N_{m+1,d}} \overline{Y_{m+1,d}^j(\zb a)}
	\\&
	\sideset{}{'}\sum_{n=m+1-l}^{l+m+1} (-1)^{\frac{n+d-2}{2}} \frac{(n-1)!!}{(n+d-3)!!\,(d-3)!!} \sum_{i=1}^{N_{n,d}} G^{n,i,d}_{m+1,j,l,k} Y_{n,d}^{i}(\zb\omega).
	\end{align*}
	Since, by \eqref{eq:VSd},
	\begin{equation*}
	\frac{\abs{\S^{d-1}}}{\abs{\S^{d-2}}}
	= \frac{\sqrt{\pi}\, \Gamma(\frac{d-1}{2})}{\Gamma(\frac{d}{2})}
	= \frac{\sqrt{\pi}\, (\frac{d-3}{2})!}{\Gamma(\frac{d}{2})}
	= \frac{\sqrt\pi\,2^{-\frac{d-3}{2}}\,(d-3)!!}{\Gamma(\frac{d}{2})},
	\end{equation*}
	we obtain
	\begin{multline*}
	\mathcal D^\text{(odd)}\widetilde V_{m,l,k}(\zb a,\zb\omega)
	\\
	= \frac{2^{\frac{d+1}{2}}\pi^{\frac{d-1}{2}}}{\sqrt{2m+d}}\, \sum_{j=1}^{N_{m+1,d}} \overline{Y_{m+1,d}^j(\zb a)}
	\sideset{}{'}\sum_{n=m+1-l}^{l+m+1} (-1)^{\frac{n+1}{2}} \frac{(n-1)!!}{(n+d-3)!!} \sum_{i=1}^{N_{n,d}} G^{n,i,d}_{m+1,j,l,k} Y_{n,d}^{i}(\zb\omega).\qedhere
	\end{multline*} 
\end{proof}

\begin{Theorem}
	\label{thm:cone-svd-norm}
	
	The functions 
	$$\frac{1}{\lambda_{m,l,d}}\, \mathcal D^\text{(odd)} \widetilde V_{m,l,k}
	,\qquad m\in\N,\; l=0,\dots,m,\;k=1,\dots,N_{l,d},\;l+m\text{ even} $$ 
	are orthonormal in $L^2(\S^{d-1}\times\S^{d-1})$, where
	\begin{equation}
	\lambda_{m,l,d}
	= \sqrt{\frac{N_{m+1,d}\,\mu_{m,d}^2\,\abs{\S^{d-2}}}{\abs{\S^{d-1}}^2}\,
	\sideset{}{'}\sum_{n=m+1-l}^{m+1+l} \nu_{n,d}^2\, N_{n,d}
	\left< P_{m+1,d}\, P_{l,d}, P_{n,d}\right>_{w_d}}
	\label{eq:lambda_mld}
	\end{equation}
	and
	\begin{equation*}
	\left< P_{m+1,d}\, P_{l,d}, P_{n,d}\right>_{w_d}
	= \int_{-1}^1 P_{m+1,d}(t)\, P_{l,d}(t)\, P_{n,d}(t)\,(1-t^2)^{\frac{d-3}{2}} \d t.
	\end{equation*}
\end{Theorem}

\begin{proof}
	Let $m,m'\in\N$, $l=0,\dots,m$, $l'=0,\dots,m'$, $k=1,\dots,N_{l,d}$, $k'=1,\dots,N_{l',d}$ such that $m+l$ and $m'+l'$ are even.
	We have
	\begin{align*}
	&\left<\mathcal D^\text{(odd)} \widetilde V_{m,l,k}, \mathcal D^\text{(odd)} \widetilde V_{m',l',k'}\right>_{L^2(\S^{d-1}\times\S^{d-1})}
	\\
	={}& \mu_{m,d}\,\mu_{m',d}\,
	\sum_{j=1}^{N_{m+1,d}} \sum_{j'=1}^{N_{m'+1,d}}\int_{\S^{d-1}}  \overline{Y_{m+1,d}^j(\zb a)}\, Y_{m'+1,d}^{j'}(\zb a) \d\zb a  
	\\&
	\sideset{}{'}\sum_{n=m+1-l}^{m+1+l} \sideset{}{'}\sum_{n'=m'+1-l'}^{m'+1+l'} \nu_{n,d}\, \nu_{n',d}
	\sum_{i=1}^{N_{n,d}}\sum_{i'=1}^{N_{n',d}} G^{n,i,d}_{m+1,j,l,k}\, \overline{G^{n',i',d}_{m'+1,j',l',k'}}
	\int_{\S^{d-1}} Y_{n,d}^{i}(\zb\omega)\, \overline{Y_{n',d}^{i'}(\zb\omega)}\d\zb\omega.
	\end{align*}
	By the orthonormality of the spherical harmonics, we obtain
	\begin{multline*}
	\left<\mathcal D^\text{(odd)} \widetilde V_{m,l,k}, \mathcal D^\text{(odd)} \widetilde V_{m',l',k'}\right>_{L^2(\S^{d-1}\times\S^{d-1})}
	\\
	= \delta_{m,m'}\, \mu_{m,d}^2\,	\sum_{j=1}^{N_{m+1,d}}
	\sideset{}{'}\sum_{n=m+1-l}^{m+1+l} \nu_{n,d}^2
	\sum_{i=1}^{N_{n,d}} G^{n,i,d}_{m+1,j,l,k}\, \overline{G^{n,i,d}_{m+1,j,l',k'}}.
	\end{multline*} 
	We have by the definition of the Gaunt coefficients in \eqref{eq:Gaunt}
	\begin{align*}
	& \sum_{j=1}^{N_{m+1,d}} \sum_{i=1}^{N_{n,d}} G^{n,i,d}_{m+1,j,l,k} \overline{G^{n,i,d}_{m+1,j,l',k'}}
	\\={}& \sum_{j=1}^{N_{m+1,d}} \sum_{i=1}^{N_{n,d}}
	\int_{\S^{d-1}} Y_{m+1,d}^j(\zb\xi)\, Y_{l,d}^k(\zb\xi)\, \overline{Y_{n,d}^i(\zb\xi)} \d\zb\xi \int_{\S^{d-1}} \overline{Y_{m+1,d}^j(\zb\eta)\, Y_{l',d}^{k'}(\zb\eta)}\, Y_{n,d}^i(\zb\eta) \d\zb\eta
	\\={}&  \frac{N_{m+1,d}\,N_{n,d}}{\abs{\S^{d-1}}^2}
	\int_{\S^{d-1}} \int_{\S^{d-1}} P_{m+1,d}(\zb\xi^\top\zb\eta)\, P_{n,d}(\zb\xi^\top\zb\eta)\, Y_{l,d}^k(\zb\xi) \d\zb\xi\, \overline{Y_{l',d}^{k'}(\zb\eta)} \d\zb\eta,
	\end{align*}
	where the last equality follows from the addition formula \eqref{eq:addition-theorem} for spherical harmonics.
	Applying the Funk--Hecke formula \eqref{eq:Funk-Hecke} to the inner integral, we obtain
	\begin{align*}
	&\sum_{j=1}^{N_{m+1,d}} \sum_{i=1}^{N_{n,d}} G^{n,i,d}_{m+1,j,l,k} \overline{G^{n,i,d}_{m+1,j,l',k'}}
	\\={}& \frac{N_{m+1,d}\,N_{n,d} \abs{\S^{d-2}}}{\abs{\S^{d-1}}^2} 
	\int_{-1}^1 P_{m+1,d}(t)\, P_{n,d}(t)\ P_{l,d}(t)\,(1-t^2)^{\frac{d-3}{2}} \d t 
	\int_{\S^{d-1}} Y_{l,d}^k(\zb\eta)\, Y_{l',d}^k(\zb\eta) \d\zb\eta
	\\={}& \delta_{l,l'}\,\delta_{k,k'}\,\frac{N_{m+1,d}\,N_{n,d}}{\abs{\S^{d-1}}^2} \abs{\S^{d-2}}
	\int_{-1}^1 P_{m+1,d}(t)\, P_{n,d}(t)\ P_{l,d}(t)\,(1-t^2)^{\frac{d-3}{2}} \d t,
	\end{align*}
	where we used again the orthonormality of the spherical harmonics.
	By \cite{ChKo10}, the value of the integral
	$
	\left< P_{m+1,d}\, P_{n,d}, P_{l,d}\right>_{w_d}
	$
	is nonzero if and only if 
	$$n \in\{ \abs{m+1-l},\, \abs{m+1-l}+2,\,\dots,\, m+1+l\}.$$
	Hence, we have
	\begin{multline*}
	\left<\mathcal D^\text{(odd)} \widetilde V_{m,l,k}, \mathcal D^\text{(odd)} \widetilde V_{m',l',k'}\right>_{L^2(\S^{d-1}\times\S^{d-1})}
	\\
	={} \delta_{m,m'}\, \delta_{l,l'}\, \delta_{k,k'}\, \frac{N_{m+1,d}\,\mu_{m,d}^2\,\abs{\S^{d-2}}}{\abs{\S^{d-1}}^2}\,
	\sideset{}{'}\sum_{n={m+1-l}}^{m+1+l} \nu_{n,d}^2\, N_{n,d}
	\left< P_{m+1,d}\, P_{n,d}, P_{l,d}\right>_{w_d}.
	\end{multline*}
\end{proof}

\subsection{Bounds on the singular values}
\label{sec:cone-beam-sv-bounds}
\subsubsection{Upper bound}
In this section, we show that the singular values
$
\lambda_{m,l,d}
$
of the cone-beam transform $\mathcal D^{(odd)}$,
which are given in \eqref{eq:lambda_mld},
are bounded independently of $m$ and $l$, which implies that the cone-beam transform as operator $\mathcal D^\text{(odd)}\colon L^2(\B^d)\to L^2(\S^{d-1}\times\S^{d-1})$ is bounded.

\begin{Lemma}
	\label{lem:nu*N}
	Let $n\ge1$ and $d\ge3$ be odd integers. Then
	\begin{equation*}
	\nu_{n,d}^2 N_{n,d} \le
	\begin{cases}
	\pi,&d=3\\
	\frac{d}{((d-2)!!)^2},& d\ge5.
	\end{cases}
	\end{equation*}
\end{Lemma}
\begin{proof}
	We have by \eqref{eq:nu} and \eqref{eq:Nnd}
	\begin{align*}
	\nu_{n,d}^2 N_{n,d}
	&= \frac{(n-1)!!^2}{(n+d-3)!!^2} \frac{(2n+d-2)\,(n+d-3)!}{n!\,(d-2)!}
	\\&= \frac{(n-1)!!}{n!!} \frac{(n+d-4)!!}{(n+d-3)!!} \frac{2n+d-2}{(d-2)!}
	.
	\end{align*}
	We note that, by \cite{Bau07}, for $n\to\infty$
	$$
	(2n-1)!! \simeq \frac{(2n)!!}{\sqrt{\pi(n+\frac12)}}.
	$$
	Hence, we obtain
	\begin{equation}
	\nu_{n,d}^2 N_{n,d}
	\simeq \frac{\sqrt{\pi(\frac{n}{2}+1)}}{n+1} \frac{\sqrt{\pi(\frac{n+d-1}{2})}}{n+d-4} \frac{2n+d-2}{(d-2)!}
	\simeq \frac{\pi}{(d-2)!}
	.
	\label{eq:nu*n-limit}
	\end{equation}
	Because
	\begin{align*}
	\frac{\nu_{n+2,d}^2 N_{n+2,d}}{\nu_{n,d}^2 N_{n,d}}
	&= \frac{n+1}{n+2} \frac{n+d-2}{n+d-1} \frac{2n+d+2}{2n+d-2}
	\\&= \frac{2n^3+3dn^2+(d^2+3d-6)n+d^2-4 \hfill{}}{2n^3+3dn^2+(d^2+3d-6)n+2d^2-6d+4}
	,
	\end{align*}
	the sequence $n\mapsto\nu_{n,d}^2 N_{n,d}$ is increasing for $d=3$ and decreasing for $d\ge5$.
	The fact that
	\begin{equation*}
	\nu_{1,d}^2 N_{1,d}
	= \frac{d}{((d-2)!!)^2}
	\end{equation*}
	completes the proof.
\end{proof}

\begin{Theorem}
	\label{thm:sv-upper-bound}
	Let $d\ge3$ be an odd integer and $m,l\in\N$ such that $l\le m$ and $m+l$ is even.
	Then the singular values $\lambda_{m,l,d}$ satisfy
	\begin{equation*}
	\lambda_{m,l,d}
	\le  {2^{\frac{d+1}{4}}\,\pi^{\frac{d-1}{4}}\, \sqrt{C_d\, (d-2)!!}\, \sqrt{\frac{l+1}{2m+d}}}
	\le (2\pi)^{\frac{d-1}{4}}\, \sqrt{C_d\, (d-2)!!}
	,
	\end{equation*}
	where 
	\begin{equation*}
	C_d = 
	\begin{cases}
	\pi,&d=3\\
	\frac{d}{((d-2)!!)^2},& d\ge5.
	\end{cases}
	\end{equation*}
	In particular, we have
	$
	\lim_{m\to\infty} \lambda_{m,l,d} = 0
	$
	for all $l\in\N$.
\end{Theorem}
\begin{proof}
Because of the orthogonality \eqref{eq:Legendre-orthogonal} of the Legendre polynomials, we have
\begin{equation*}
P_{l,d}\, P_{n,d}
= \sum_{m=\abs{l-n}-1}^{l+n-1} \frac{N_{m+1,d}\abs{\S^{d-2}}}{\abs{\S^{d-1}}} \left<P_{m+1,d}P_{l,d}, P_{n,d} \right>_{w_d} P_{m+1,d}.
\end{equation*}
Utilizing the fact that $P_{i,d}(1)=1$ for all $i\in\N$, we obtain
\begin{equation}
1
= \sum_{m=\abs{l-n}-1}^{l+n-1} \frac{N_{m+1,d}\abs{\S^{d-2}}}{\abs{\S^{d-1}}} \left<P_{m+1,d}P_{l,d}, P_{n,d} \right>_{w_d}.
\label{eq:sum-1}
\end{equation}
Since all summands in the above sum \eqref{eq:sum-1} are non-negative, they are bounded by
\begin{equation}
\frac{N_{m+1,d}\abs{\S^{d-2}}}{\abs{\S^{d-1}}} \left<P_{m+1}P_l, P_n \right>_{w_d} \le1.
\label{eq:summands<1}
\end{equation}
Inserting the bound from \prettyref{lem:nu*N} into the definition of the singular values \eqref{eq:lambda_mld}, we have
\begin{align*}
\lambda_{m,l,d}^2
&= \frac{N_{m+1,d}\,\mu_{m,d}^2\,\abs{\S^{d-2}}}{\abs{\S^{d-1}}^2}\,
\sideset{}{'}\sum_{n=m+1-l}^{m+1+l} \nu_{n,d}^2\, N_{n,d}
\left< P_{m+1,d}\, P_{l,d}, P_{n,d} \right>_{w_d}
\\
&\le \frac{N_{m+1,d}\,\mu_{m,d}^2\,\abs{\S^{d-2}}}{\abs{\S^{d-1}}^2}\, C_d
\sideset{}{'}\sum_{n=m+1-l}^{m+1+l} 
\left< P_{m+1,d}\, P_{l,d}, P_{n,d} \right>_{w_d}.
\end{align*}
With \eqref{eq:summands<1}, we obtain
\begin{align*}
\lambda_{m,l,d}^2
& \le  C_d\, \frac{\mu_{m,d}^2}{\abs{\S^{d-1}}}\,
\sideset{}{'}\sum_{n=m+1-l}^{m+1+l} 1
\\& = C_d\, \frac{\mu_{m,d}^2}{\abs{\S^{d-1}}} \, (l+1)
\\& = C_d\, {2^{\frac{d+1}{2}}\,\pi^{\frac{d-1}{2}}\,(d-2)!!}\, \frac{l+1}{2m+d}
,
\end{align*}
where we inserted the formulas of $\mu_{m,d}$ from \eqref{eq:mu} and $\abs{\S^{d-1}}$ from \eqref{eq:VSd}.
\end{proof}

\subsubsection{Lower bound}
\begin{Theorem}
	\label{thm:sv-lower-bound}
	Let $d\ge3$ be an odd integer.
	There exists a constant $c_d>0$, which depends only on the dimension $d$, such that for all $m\in\N$ and $l\in\{0,\dots,m\}$ with $m+l$ even, the singular values admit the lower bound
	\begin{equation}
	\abs{\lambda_{m,l,d}}
	\ge c_d\, m^{-1/2}.
	\label{eq:sv-lower-bound}
	\end{equation}
	This bound is asymptotically tight, in the sense that the exponent $-1/2$ in \eqref{eq:sv-lower-bound} cannot be replaced by a greater one.
\end{Theorem}
\begin{proof}
We extract the smallest value of $\nu_{n,d}^2$ in the following sum
\begin{align*}
\lambda_{m,l,d}^2
&= \frac{N_{m+1,d}\,\mu_{m,d}^2\,\abs{\S^{d-2}}}{\abs{\S^{d-1}}^2}\,
\sideset{}{'}\sum_{n=m+1-l}^{m+1+l} \nu_{n,d}^2\, N_{n,d}
\left< P_{m+1,d}\, P_{l,d}, P_{n,d} \right>_{w_d}
\\
&\ge \frac{N_{m+1,d}\,\mu_{m,d}^2}{\abs{\S^{d-1}}}\,
\left(\min_{n=m+1-l}^{m+1+l}\nu_{n,d}^2\right) \sideset{}{'}\sum_{n=m+1-l}^{m+1+l} \, \frac{N_{n,d}\, \abs{\S^{d-2}}}{\abs{\S^{d-1}}}
\left< P_{m+1,d}\, P_{l,d}, P_{n,d} \right>_{w_d}.
\end{align*}
Utilizing \eqref{eq:sum-1} with the roles of $m+1$ and $n$ interchanged, we obtain
\begin{equation*}
\lambda_{m,l,d}^2
\ge \frac{N_{m+1,d}\,\mu_{m,d}^2}{\abs{\S^{d-1}}}\,
\min_{n=m+1-l}^{m+1+l}\nu_{n,d}^2.
\end{equation*}
Since the map
\begin{equation*}
n\mapsto
\nu_{n,d}^2
= \frac{(n-1)!!^2}{(n+d-3)!!^2}
\end{equation*}
is decreasing, we have
\begin{equation*}
\min_{n=m+1-l}^{m+1+l}\nu_{n,d}^2
= \nu_{m+1+l,d}^2
.
\end{equation*}
Because $0\le l \le m$ and again $\nu_{m+1+l,d}^2$ decreases with respect to $l$, we further see that
\begin{equation*}
\min_{n=m+1-l}^{m+1+l}\nu_{n,d}^2
\ge \nu_{2m+1,d}^2
= \frac{(2m)!!^2}{(2m+d-2)!!^2}.
\end{equation*}
Hence, we have
\begin{align*}
\lambda_{m,l,d}^2
&\ge \frac{N_{m+1,d}\,\mu_{m,d}^2}{\abs{\S^{d-1}}}\,
\frac{(2m)!!^2}{(2m+d-2)!!^2}
\\
&= \frac{2^{d+1}\,\pi^{d-1}}{(2m+d)}\,
\frac{(d-2)!!}{2^{\frac{d+1}{2}} \pi^{\frac{d-1}{2}}}\,
\frac{(2m+d)\,(m+d-2)!}{(m+1)!\,(d-2)!}\,
\frac{(2m)!!^2}{(2m+d-2)!!^2}
\\
&= \frac{2^{\frac{d+1}{2}} \pi^{\frac{d-1}{2}}}{(d-3)!!}\,
\frac{(m+d-2)!}{(m+1)!}\,
\frac{(2m)!!^2}{(2m+d-2)!!^2},
\end{align*}
where we inserted \eqref{eq:Nnd}, \eqref{eq:mu} and \eqref{eq:VSd}.
We are going to apply Stirling's approximation of the factorial
\begin{equation*}
n! \simeq \sqrt{2\pi}\,n^{n+1/2}\,\e^{-n}
\end{equation*}
and the double factorials, cf.\ \cite{Bau07},
\begin{align*}
(2n)!! &\simeq \sqrt{\pi}\, ({2n})^{n+1/2} \,\e^{-n} ,
&
(2n-1)!! &\simeq \sqrt{2}\, ({2n})^n\,\e^{-n}.
\end{align*}
We obtain for $m\to\infty$
\begin{align*}
\lambda_{m,l,d}^2
&\ge\frac{2^{\frac{d+1}{2}} \pi^{\frac{d-1}{2}}}{(d-3)!!}\,
\frac{(m+d-2)!}{(m+1)!}\,
\frac{(2m)!!^2}{(2m+d-2)!!^2}
\\&\simeq
\frac{2^{\frac{d+1}{2}} \pi^{\frac{d-1}{2}}}{(d-3)!!}\,
\frac{(m+d-2)^{m+d-3/2}\,\e^{-m-d+2}}{(m+1)^{m+3/2}\,\e^{-m-1}}\,
\frac{(2m)^{2m+1}\,\e^{-2m} \pi}{2(2m+d-1)^{2m+d-1}\,\e^{-2m-d-1}}
\\&\simeq
\frac{2^{\frac{3-d}{2}} \pi^{\frac{d+1}{2}}}{(d-3)!!}\, \e^{4}
\frac{(m+d-2)^{m+d-3/2}}{(m+1)^{m+3/2}}\,
\frac{m^{2m+1}}{(m+\frac{d-1}{2})^{2m+d-1}}.
\end{align*}
Hence, there exists a constant $c_d>0$ such that 
\begin{equation*}
\lambda_{m,l,d}
\ge
\sqrt{\frac{2^{\frac{d+1}{2}} \pi^{\frac{d-1}{2}}}{(d-3)!!}\,
\frac{(m+d-2)!}{(m+1)!}\,
\frac{(2m)!!^2}{(2m+d-2)!!^2}}
\ge c_d\, m^{-1/2}.
\end{equation*}
In order to show that this bound is tight, we consider the case $m$ even and $l=0$. We have by \eqref{eq:Legendre-orthogonal}
\begin{align*}
\lambda_{m,0,d}^2
&= \frac{N_{m+1,d}\,\mu_{m,d}^2\,\abs{\S^{d-2}}}{\abs{\S^{d-1}}^2}\,
\nu_{m+1,d}^2\, N_{m+1,d}
\left< P_{m+1,d}, P_{m+1,d} \right>_{w_d}
\\
&= \frac{N_{m+1,d}\,\mu_{m,d}^2}{\abs{\S^{d-1}}}\,
\nu_{m+1,d}^2
.
\end{align*}
By \prettyref{eq:nu*n-limit}, we have for $m\to\infty$
\begin{align*}
\lambda_{m,0,d}^2
& \simeq \frac{\pi}{(d-2)!} \frac{\mu_{m,d}^2}{\abs{\S^{d-1}}}.\qedhere
\end{align*}
\end{proof}

\begin{Remark}
	While the lower bound $\bigo(m^{-1/2})$ 
	on the singular values $\lambda_{m,l,d}$ is asymptotically strict, we have only shown that they can be bounded from above by a constant in \prettyref{thm:sv-upper-bound}.
	However, the degree of ill-posedness of the reconstruction problem depends on the behavior of the smallest singular values, which is here $\bigo(m^{-1/2})$ and so the same as for the Radon transform in 2D.
\end{Remark}

\section{Cone-beam transform in $\R^3$}
\label{sec:R3}

In this section, we state the singular value decomposition of the cone-beam transform $\mathcal D^{(odd)}$ from \prettyref{sec:cone-beam} for the dimension $d=3$.
This case was already shown in \cite{Kaz15}.
Before we state the result, we give some formulas of spherical harmonics and Gaunt coefficients in this case.
We write a point $\zb\xi\in\S^2$ in cylindrical coordinates
\begin{equation*}
\zb\xi(\varphi,t)=(\cos\varphi \,\sqrt{1-t^2},\sin\varphi \,\sqrt{1-t^2},t)^{\top}, \quad\varphi\in[0,2\pi),\, t\in[-1,1].
\end{equation*}
We define the normalized associated Legendre functions of degree $n\in\N$ and order $k=-n,\dots,n$ by
\begin{equation*}
\widetilde P_n^k
= \sqrt{\frac{2n+1}{4\pi}\,\frac{(n-k)!}{(n+k)!}}\, 
\frac{(-1)^{k}}{2^{n}n!}\left(1-t^{2}\right)^{k/2}\frac{\mathrm{d}^{n+k}}{\mathrm{d}t^{n+k}}\left(t^{2}-1\right)^{n}
.
\end{equation*}
The spherical harmonics
\begin{equation*}
Y_{n}^{k}(\zb\xi(\varphi,t))=
\widetilde P_n^k(t)\, \mathrm{e}^{\mathrm{i}k\varphi}
,\qquad \zb \xi(\varphi,t)\in\S^2,
\end{equation*}
of degree $n\in\N$ and order $k\in\{-n,\dots,n\}$
form an orthonormal basis of $L^{2}(\S^2)$, see \cite[Chapter 5]{Varsha88} and also \cite{DaXu2013}.
The Gaunt coefficients
\begin{equation}
G^{n,k}_{n_1,k_1,n_2,k_2}
= \int_{\S^2} Y_{n_1}^{k_1}(\zb\xi)\, Y_{n_2}^{k_2}(\zb\xi)\, \overline{Y_{n}^{k}(\zb\xi)} \d\zb\xi
\label{eq:Gaunt}
\end{equation}
are zero unless all the conditions
$$k= k_1+k_2,\; \abs{k_1}\le n_1,\; \abs{k_2}\le n_2,\; \abs{k}\le n$$
and
$$n=\abs{n_1-n_2},\abs{n_1-n_2}+2,\dots, n_1+n_2$$
hold.
An explicit representation of Gaunt coefficients can be found in \cite{Gau29}.
The Gaunt coefficients are closely related to the Clebsch--Gordan coefficients, cf.\ \cite{Varsha88}.

\begin{Theorem}
	The odd cone-beam transform $\mathcal D^{(odd)}\colon L^2(\B^3)\to L^2(\S^{2}\times\S^{2})$ has the singular value decomposition 
	$$
	\mathcal D^{(\mathrm{odd})} \widetilde V_{m,l,k} = \lambda_{m,l}\,W_{m,l,k}
	,\qquad m\in\N,\; 0\le l\le m,\; \text{$m+l$ even},\; k\in\{-l,\dots,l\}.
	$$
	The polynomials
	\begin{equation*}
	\widetilde V_{m,l,k}(s\zb\omega)
	= \sqrt{2m+3}\, s^l P_{\frac{m-l}{2}}^{\left(0, l+\frac{1}{2}\right)} (2s^2-1)\, Y_{l}^k(\zb\omega)
	,\qquad s\in[0,1],\; \zb\omega\in\S^{2},
	\end{equation*}
	form an orthonormal basis of $L^2(\B^3)$.
	The singular values are given by
	\begin{equation*}
	\lambda_{m,l}
	= \sqrt{2\pi\,
	\sideset{}{'}\sum_{n=m+1-l}^{m+1+l} \frac{(2n+1)(n-1)!!^2}{n!!^2}
	\left< P_{m+1}\, P_{n}, P_{l}\right>}
	,
	\end{equation*}
	where $\sum'$ denotes the summation over odd indices and
	\begin{equation*}
	\left< P_{m+1}\, P_{n}, P_{l}\right>
	=\frac{2 (l+m-n)!!\, (l-m+n-2)!!\, (-l+m+n)!!\, (l+m+n+1)!!}{(l+m-n+1)!!\, (l-m+n-1)!!\, (-l+m+n+1)!!\, (l+m+n+2)!!}
	\end{equation*}
	for $n \in\{ \abs{m+1-l},\, \abs{m+1-l}+2,\,\dots,\, m+1+l\}$ and zero otherwise.
	The singular values satisfy
	\begin{equation}
	c_1 m^{-1/2}
	\le \lambda_{m,l}
	\le c_2 m^{-1/8}
	\label{eq:lambda_ml-bound}
	\end{equation}
	for some constants $c_1$, $c_2$ independent of $m$ and $l$.
	Furthermore, the functions
	\begin{multline*}
	W_{m,l,k}(\zb a,\zb\omega) \\= \frac{4\pi}{\lambda_{m,l}\,\sqrt{2m+3}}\, \sum_{j=-m-1}^{m+1} \overline{Y_{m+1}^j(\zb a)}
	\sideset{}{'}\sum_{n=m+1-l}^{m+1+l}  \frac{(-1)^{\frac{n+1}{2}}\,(n-1)!!}{n!!} G^{n,j+k}_{m+1,j,l,k} Y_{n}^{j+k}(\zb\omega)
	\end{multline*}
	for $\zb a,\zb\omega\in\S^{2}$
	are orthonormal in $L^2(\S^2\times\S^2)$.
\end{Theorem}

\begin{proof}
	The singular value decomposition is a special case of Theorems \ref{thm:cone-svd} and \ref{thm:cone-svd-norm}.
	The triple product $\left< P_{m+1}\, P_{n}, P_{l}\right>$ is computed in \cite{Neu1878} (see also \cite{Sal56}).
	The lower bound of the singular values $\lambda_{m,l}$ in \eqref{eq:lambda_ml-bound} is due to \prettyref{thm:sv-lower-bound}.
	It is left to show the upper bound in \eqref{eq:lambda_ml-bound}, which we do as in \cite{Kaz15}.
	Changing roles of $m$ and $n$ in \eqref{eq:sum-1}, we have
	\begin{equation}
	1
	= \sum_{n=m+1-l}^{m+1+l} \frac{2n+1}{2} \left<P_{m+1}P_l, P_n \right>.
	\label{eq:sum-1-m}
	\end{equation}
	Furthermore, since $\abs{P_n(t)}\le1$ for $\abs t\le1$, we obtain the inequality
	\begin{equation*}
	\left<P_{m+1}P_l, P_n \right>
	\le \int_{-1}^1 \abs{P_{m+1,d}(t)\, P_{n,d}(t)\ P_{l,d}(t)} \d t
	\le \int_{-1}^1 \abs{P_{m+1,d}(t)\, P_{n,d}(t)} \d t
	.
	\end{equation*}
	By the Cauchy-Schwarz inequality, 
	\begin{equation}
	\left<P_{m+1}P_l, P_n \right>
	\le \norm{P_{m+1}}_{L^2(-1,1)} \norm{P_{n}}_{L^2(-1,1)}
	= \frac{2}{\sqrt{(2m+3)\, (2n+1)}}.
	\label{eq:LP-upper-bound}
	\end{equation}
	Based on Wallis' product, it was shown in \cite{HiPoQu18} that for $m\in\N$
	\begin{equation*}
	\frac{(2m)!!^2}{(2m+1)!!^2} \le \frac{\pi}{2(2m+1)}.
	\end{equation*}
	Hence, we have
	\begin{align*}
	\lambda_{m,l}^2
	&= {2\pi\, \sideset{}{'}\sum_{n=1}^{m+1+l} \frac{(2n+1)(n-1)!!^2}{n!!^2}
		\left< P_{m+1}\, P_{n}, P_{l}\right>}
	\\&\le 2\pi^2\, \sideset{}{'}\sum_{n=1}^{2m+1} \frac{2n+1}{2n}
	\left< P_{m+1}\, P_{n}, P_{l}\right>
	\end{align*}
	because $0\le l\le m$.
	By the Cauchy-Schwarz inequality, we have 
	\begin{equation*}
	\lambda_{m,l}^2
	\le 2\pi^2\, 
	\sqrt{\sideset{}{'}\sum_{n=1}^{2m+1} \frac{2n+1}{2} \left< P_{m+1}\, P_{n}, P_{l}\right>}
	\sqrt{\sideset{}{'}\sum_{n=1}^{2m+1} \frac{2n+1}{2n^2} \left< P_{m+1}\, P_{n}, P_{l}\right>}
	\end{equation*}
	Inserting \eqref{eq:sum-1-m} and \eqref{eq:LP-upper-bound}, we obtain
	\begin{align*}
	\lambda_{m,l}^2
	&= 2\pi^2\,
	\sqrt{\sideset{}{'}\sum_{n=1}^{2m+1} \frac{2n+1}{2n^2} \left< P_{m+1}\, P_{n}, P_{l}\right>}
	\\
	&\le \frac{2\pi^2}{(2m+3)^{\frac14}}\,
	\sqrt{\sideset{}{'}\sum_{n=1}^{2m+1} \frac{\sqrt{2n+1}}{n^2} }.
	\end{align*}
	The last sum converges for $m\to\infty$  and thus can be bounded from above by a constant intependent of $m$, which implies that $\lambda_{m,l}^2\in\bigo(m^{-1/4})$.
\end{proof}

\begin{Remark}
	The upper bound on the singular values $\lambda_{m,l} \in\bigo(m^{-1/8})$ may not be optimal.
	The authors think that the upper bound can be improved to $\bigo(m^{-1/2})$ even in general dimension $d$.
	This conjecture is backed by numerical computations as well as the following observation, which is not a proof though.
	We consider Grangeat's formula \eqref{eq:Grangeat}
	\begin{equation*}
	(-1)^{d}\left(\frac{\partial}{\partial s}\right)^{d-2} \mathcal Rf(\zb\omega, \zb a^\top\zb\omega)
	= \mathcal S^{(d-2)}_d \mathcal Df(\zb a,\zb\omega).
	\end{equation*}
	We know that the singular values of the Radon transform $\mathcal R$ are $\bigo(m^{(1-d)/2})$ and those of the $d-2$ differentiations are $\bigo(m^{d-2})$, so the left side should behave like $\bigo(m^{(d-3)/2})$.
	On the right side, $\mathcal S^{(d-2)}_d$ has the singular values $\bigo(m^{(d-2)/2})$ by \prettyref{lem:Sj-ev-d-approx}, so the singular values of the cone-beam transform $\mathcal D$ should behave like $\bigo(m^{-1/2})$.
\end{Remark}

\end{document}